\documentclass[leqno,10pt,twoside,a4,bezier]{article}
\usepackage{amssymb,amsmath,amsfonts,latexsym,ifthen,makeidx,bezier,epic,eepic}
\usepackage{diagrams}
\diagramstyle[tight,centredisplay%
,dpi=600
,noPostScript%
,heads=LaTeX%
]%
\newcommand{\kla}[1]{ {\langle #1 \rangle} }

\newcommand{\st}{\;|\;}
\newcommand{\mdf}{\stackrel{\rm def}{=}}

\newcommand{\dom}{ {\rm dom} }
\newcommand{\ran}{ {\rm ran} }

\newcommand{\height}{ {\rm ht} }

\newcommand{\rsubset}{\subsetneqq}
\newcommand{\sub}{\subseteq}

\newfont{\ssi}{cmssi12 at 12pt}
\newcommand{\qO}{{\Bbb Q}}

\newcommand{\rest}{\restriction}

\newcommand{\verl}{{{}^\frown}}
\newcommand{\leer}{\emptyset}
\newcommand{\ohne}{\setminus}

\newcommand{\qed}{\nopagebreak\hspace*{\fill}$ \Box $\\}

\newenvironment{ea*}{\begin{eqnarray*}}{\end{eqnarray*}}
\newtheorem{thm}{Theorem}[section]

\newtheorem{lem}[thm]{Lemma}
\newtheorem{definition}[thm]{Definition}
\newtheorem{obs}[thm]{Observation}
\newtheorem{question}[thm]{Question}
\newcommand{\claim}[2]{
     \begin{enumerate}
       \item[{#1}] {\em #2}
     \end{enumerate}}

\setbox0=\hbox{$\longrightarrow$}
\newcommand{\To}{\longrightarrow}

\newcommand{\proof}{\noindent{\em Proof.}\hspace{1em}}

\newcommand{\balpha}{{\bar{\alpha}}}

\newcommand{\tp}{{\tilde{p}}}

\newcommand{\bp}{{\bar{p}}}
\newcommand{\bq}{{\bar{q}}}

\newcommand{\vb}{{\vec{b}}}

\newcommand{\vp}{{\vec{p}}}
\newcommand{\vq}{{\vec{q}}}
\newcommand{\vs}{{\vec{s}}}
\newcommand{\vt}{{\vec{t}}}

\newcommand{\vr}{{\vec{r}}}

\newcommand{\vgamma}{{\vec{\gamma}}}

\renewcommand{\phi}{\varphi}

\newcommand{\tx}[1]{\;\mbox{\rm #1}\;}

\newcommand{\ZFC}{{\tt ZFC}}

\newcommand{\forces}{\Vdash}
\newcommand{\incomp}{\perp}
\newcommand{\Aut}{\mathop{\hbox{Aut}}}
\def\<#1>{\langle#1\rangle}
\newcommand{\restrict}{\upharpoonright}
\renewcommand{\P}{\mathord{\Bbb P}}
\newcommand{\of}{\subseteq}
\newcommand{\intersect}{\cap}
\newcommand{\union}{\cup}
\newcommand{\Union}{\bigcup}
\newcommand{\cross}{\times}
\newcommand{\compose}{\circ}

\begin{document}
\hyphenation{covering}

%
\title{Degrees of Rigidity for Souslin Trees}
\author{Gunter Fuchs\\
  Westf\"{a}lische Wilhelms-Universit\"{a}t M\"{u}nster, Germany\\
  \and
  Joel David Hamkins\thanks{The second author is grateful for the generous support of the DFG and the Westf\"alische Wilhelms-Universit\"at
M\"unster for supporting him as Mercator Gastprofessor in the summer 2004 and for the support of NWO Bezoekersbeurs B 62-612 at Universiteit van
Amsterdam in Summer 2005. In addition, he has been partially supported by
PSC-CUNY research grants. He is affiliated with the College of Staten Island of CUNY and The CUNY Graduate Center.}\\
  The City University of New York, USA}
\date{\today}
\maketitle
\begin{abstract}
We investigate various strong notions of rigidity for Souslin trees, separating them under $\diamondsuit$ into a hierarchy. Applying our methods to
the automorphism tower problem in group theory, we show under $\diamondsuit$ that there is a group whose automorphism tower is highly malleable by
forcing.
\end{abstract}

\section{Introduction}

Automorphisms and isomorphisms of $\omega_1$-trees have been long studied in set theory (see \cite{TSP}, \cite{AutomorphismsOfTrees},
\cite{ForcingWithTrees},\cite{CRAT}, \cite{ITAT},\cite{GaifmanSpecker1964}), and several of these authors have investigated various strong forms of
rigidity of such trees. Here, by considering the absoluteness of rigidity and strong rigidity properties of a tree to various forcing extensions, we
introduce several new rigidity concepts---all of which are exhibited by the generic Souslin trees added by the usual forcing---and separate them
under $\diamondsuit$ into a proper implication hierarchy.

Our motivation for looking at these particular rigidity properties arose in connection with the automorphism tower problem in group theory.
Specifically, the main result of \cite{HamkinsThomas2000:ChangingHeights} had made essential use of the absolute rigidity properties of generic
Souslin trees to construct in the corresponding Souslin tree forcing extension a group whose automorphism tower is highly malleable by forcing. In
the final section of this paper, we replace the Souslin tree forcing argument of \cite{HamkinsThomas2000:ChangingHeights} with a construction from
$\diamondsuit$, and we conclude, consequently, that there are such groups in the constructible universe $L$.

Before introducing the rigidity notions in which we are interested, let us first set some background terminology, which we hope most readers will
find familiar. Specifically, a {\it tree} is a partial order $T=\kla{T,<_T}$ in which the predecessors of any node are well ordered and there is a
unique minimal element called the {\it root}. We will usually conflate the tree with its underlying set. The {\it height} of a node $p$ in $T$,
denoted $|p|$ or $|p|_T$, is the order type of its predecessors. We write $T(\alpha)$ for the {\it $\alpha^{\rm th}$ level} of $T$, the set of nodes
having height $\alpha$. The {\it height} of a tree $T$, $\height(T)$, is the supremum of the successors of the heights of its nodes.
%
%
We write $T|\alpha$ for the subtree of $T$ of nodes having height less than $\alpha$, and more generally, for any set $S$ of ordinals, $T|S$ is the
suborder of $T$ consisting of the nodes on a level in $S$. For $\alpha\leq\omega_1$, an $\alpha$-tree is a tree of height $\alpha$ with all levels
countable. Such a tree is {\it normal} if every node has (at least) two immediate successors (except those on the top level, if $\alpha$ is a
successor ordinal), nodes on limit levels are uniquely determined by their sets of predecessors, and every node has successors on all higher levels
up to $\alpha$. A tree is {\it $\gamma$-splitting} if every node has exactly $\gamma$ many immediate successors. It is {\it uniform} if it is
$\gamma$-splitting, for some $\gamma$. We write $T_p$ to denote the subtree of $T$ consisting of the nodes $q\in T$ with $q\ge_T p$. We write
$T_{[p]}$ for the subtree of $T$ consisting of the nodes $q\in T$ that are comparable with $p$. A branch in $T$ is a maximal linearly ordered subset
of $T$, and the {\it length} of the branch is its order type. We write $[T]$ for the set of cofinal branches, those branches containing nodes on
every level. A tree $T$ is {\it Aronszajn} if it is a normal $\omega_1$-tree with no cofinal branch. An {\it antichain} in a tree is a set of
pairwise incomparable elements. A {\it Souslin} tree is a normal $\omega_1$-tree with no uncountable antichain. When forcing with a tree, we reverse
the order, so that stronger conditions are higher up in the tree. Consequently, Souslin trees are c.c.c.~as notions of forcing. It is well known that
they are also countably distributive (see \cite[Lemma 15.28]{Jech:SetTheory3rdEdition}). An {\it automorphism} of a tree is an isomorphism of the
tree with itself; it is {\it nontrivial} if it is not the identity function. We are now ready to define the various rigidity notions.

\begin{definition}\rm
Suppose that $T$ is a tree.
\begin{enumerate}
  \item $T$ is {\it rigid} if there is no nontrivial automorphism of $T$.
  \item $T$ is {\it totally rigid} if whenever $p$ and $q$ are distinct nodes in $T$, then $T_p$ and $T_q$ are not isomorphic.
  \item $T$ has the {\it unique branch property} (UBP) if $1\forces_T T$
    has exactly one new cofinal branch.
  \item $T$ is {\it absolutely rigid} if $1\forces_T T$ is rigid.
  \item $T$ is {\it absolutely totally rigid} if $1\forces_T T$ is totally rigid.
  \item $T$ is {\it absolutely UBP} if $1\forces_T T$ has the
  UBP. Equivalently, forcing with $T\cross T$ adds precisely 2 new cofinal
  branches.\footnote{See the proof of theorem \ref{implications}.}
  \item For any property $P$, we say that $T$ is {\it absolutely $P$} if $1\forces_T T$ has property $P$, and more generally,
          $T$ is $\qO$-absolutely $P$ if $1\forces_\qO T$ has property $P$.
\end{enumerate}
\end{definition}

We shall be interested primarily in the rigidity properties of uniform normal trees, because it is (too) easy to construct rigid non-uniform trees,
simply by insisting that nodes on a level have distinct numbers of successors. But a simple back-and-forth argument shows that no uniform normal tree
of countable height can be rigid -- see \cite[Lemma 3.6.]{HamkinsThomas2000:ChangingHeights}; the argument dates back to \cite{Kurepa1935}. Abraham
\cite{CRAT} showed (in ZFC) that a rigid Aronszajn tree exists. Note that forcing with an absolutely rigid tree must preserve $\omega_1$, since
otherwise the tree would have countable height in the extension and hence fail to be rigid there. So the existence of an absolutely rigid Aronszajn
tree cannot be proved from ZFC alone, since it is consistent that every Aronszajn tree is special and thus collapses $\omega_1$. Similarly, forcing
with a normal tree with the unique branch property also preserves $\omega_1$, because every countable normal tree has continuum many branches. We
will therefore concentrate, in our $\diamondsuit$ constructions, on building Souslin trees, the canonical candidates for trees preserving $\omega_1$.

We now observe some elementary implications between these rigidity notions.

\goodbreak
\begin{thm}
\label{implications} Suppose that $T$ is a normal tree.
\begin{enumerate}
  \item If\/ $T$ is absolutely rigid or totally rigid, then it is rigid.
  \item If\/ $T$ is absolutely totally rigid then it is totally rigid and absolutely rigid.
  \item If\/ $T$ is UBP, then it is totally rigid.
  \item If\/ $T$ is absolutely UBP, then it is UBP and absolutely totally rigid.
\end{enumerate}
\end{thm}
\begin{figure}[htb]
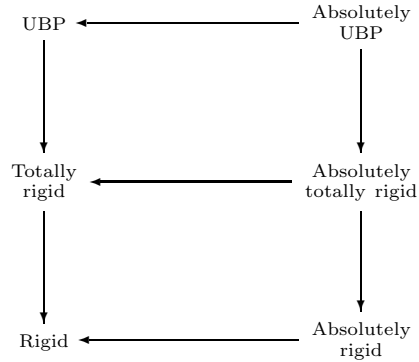

\begin{diagram}[width=6em,height=3em,objectstyle=\rm,textflow]
 {\scriptstyle UBP}             &   \lTo        &  {Absolutely\atop UBP}               \\
 \dTo            &               &  \dTo                         \\
 {Totally\atop rigid}   &   \lTo        &  {Absolutely\atop totally\ rigid}     \\
  \dTo           &               &  \dTo                         \\
 {\scriptstyle Rigid}           &  \lTo         &  {Absolutely\atop rigid}             \\
\end{diagram}
\caption{Implication Diagram} \label{fig:ImplicationDiagram}
\end{figure}
\proof If $T$ is absolutely rigid, then it is rigid, because a nontrivial automorphism of $T$ in the ground model would still be a nontrivial
automorphism in the extension. If $T$ is totally rigid, then it is rigid, because a nontrivial automorphism $\pi$ of $T$ that moves, say, $p$ to $q$,
would yield $\pi\rest T_p:T_p\cong T_q$, contradicting total rigidity. So 1 holds. Similar reasoning shows that if $T$ is absolutely totally rigid,
then it is totally rigid, and also absolutely rigid, so 2 holds.

Suppose next that $T$ has the unique branch property. If $T$ were not totally rigid, then there would be distinct nodes $p$ and $q$ with an
isomorphism $\pi:T_p\cong T_q$. We may assume without loss of generality that $|p|\leq|q|$. Thus, there is some extension $p'\geq p$ such that
$p'\perp q$. Let $b$ be $V$-generic for $T$ with $p'\in b$. Thus, $b$ provides a branch through $T_p$, and consequently $\pi[b]$ is a branch through
$T_q$. Since $p'\perp q$, these branches are distinct, and so in $V[b]$, there are at least two branches through $T$, contradicting the unique branch
property. So 3 holds.

Before continuing, let us first explain why the two definitions we gave of the absolute unique branch property are equivalent. If forcing with
$T\cross T$ necessarily adds only two new cofinal branches, then forcing with $T$ and then forcing with $T$ again clearly adds exactly one new branch
each time, and so $1\forces_T(1\forces_T T$ has exactly one new cofinal branch$)$. Conversely, it is enough to show that $T$ has the unique branch
property. But if this were not the case, then would be a name $\tau$ and a condition $p$ forcing that it is a new cofinal branch, different from the
generic branch.\footnote{When referring to {\it the} generic branch of a tree, we mean the canonical generic branch, the branch derived from the
generic filter. So here, $p\Vdash\tau\neq\Gamma$, where $\Gamma$ is the canonical name for the generic object. When $T$ is Souslin, of course, {\it
every} cofinal branch of $T$ is generic, because every antichain is bounded by a level of the tree. See \cite[Thm. II.4.7]{TSP}.} In particular, one
can strengthen $p$ so as to decide $\tau$ in various incompatible ways. If $V[b_1][b_2]$ is obtained by forcing with $T\cross T$ and $p\in b_2$, then
it follows by a simple density argument that $\tau_{b_2}$ will be a cofinal branch not in $V[b_1]$, and different from $b_2$, contradicting
$1\forces_T(1\forces_T T$ has exactly one new cofinal branch$)$. So the two definitions are in fact equivalent.

To prove 4, suppose $T$ has the absolute unique branch property. Thus, forcing with $T\cross T$ adds precisely two new branches $b_1$ and $b_2$
through $T$. If $V[b_1]$ already had 2 new cofinal branches through $T$, then there would be at least three such branches in $V[b_1,b_2]$, contrary
to our assumption, and so $T$ has the unique branch property in $V$. Lastly, we argue that $T$ is absolutely totally rigid. If not, then in some
extension $V[b]$ obtained by forcing with $T$, there would be an isomorphism $T_p\cong T_q$. As above, we may assume $p\perp q$. Further forcing to
add a $V[b]$-generic branch $c$ containing $p$ will also add an isomorphic copy of the branch through $q$, resulting that $V[b,c]$ has at least three
new cofinal branches through $T$, contradicting our assumption that $T$ was absolutely UBP. \qed

Our main result is that no other implications are provable in \ZFC.

\newtheorem{mainthm}{Main Theorem}[section]
\begin{mainthm}\label{thm:MainTheorem}
The implication diagram of Figure \ref{fig:ImplicationDiagram} is complete. Namely, if \ZFC\ is consistent, then no implication relations other than
those appearing in the transitive closure of Figure \ref{fig:ImplicationDiagram} are provable in \ZFC. Indeed, if $\diamondsuit$ holds, then there
are Souslin trees exhibiting each of the non-implications of Figure \ref{fig:ImplicationDiagram}.
\end{mainthm}

This theorem will be proved in Section \ref{Section:SeparatingTheRigidityNotions}. In Section \ref{Section:LargerContext}, we place the diagram in a
larger context including many other rigidity notions.

\section{Very rigid Souslin trees exist generically and under $\diamondsuit$}

Before proving that the implication diagram is complete, let us briefly show that the rigidity notions we have introduced actually occur. We will
show that the generic Souslin trees one adds by the usual forcing exhibit all of the rigidity notions that we have mentioned above and more. And we
will also construct such highly rigid trees under the hypothesis of $\diamondsuit$.

The usual forcing to add a Souslin tree is the partial order $\P$ consisting of all normal $\alpha$-trees, subtrees of $2^{<\alpha}$, for any
countable ordinal $\alpha$, ordered by end-extension. If $G\subset\P$ is a $V$-generic filter, then the resulting tree $T=\bigcup G$ is called a {\it
$V$-generic Souslin tree}. This forcing is countably closed. Souslin trees were first added by forcing in 1964 by Tennenbaum (see
\cite{Tennenbaum1968:SouslinsProblem}) and by Jech \cite{Jech1967:NonprovabilityOfSH}.

\begin{thm}
Every $V$-generic Souslin tree $T$ is in fact a Souslin tree in $V[T]$, and exhibits in $V[T]$ all of the rigidity properties appearing in Figure
\ref{fig:ImplicationDiagram}. \label{thm:GenericSouslinTreesAreRigid}
\end{thm}

\proof Most of this proof is very well known (see \cite[Thm. 15.23]{ST3}, \cite{TSP}); we include it for completeness and because it will motivate
some of our later constructions. Suppose that $T=\bigcup G$ for a $V$-generic filter $G\of\P$. It is easy to see that $T$ is in fact an
$\omega_1$-tree. Suppose that $A\of T$ is a maximal antichain of $T$ in $V[G]$. Let $\dot A$ be a name for $A$, such that some condition
$t_0\forces\dot A$ is a maximal antichain in the generic tree $T$. For any condition $t$ in $\P$, we may use the fact that $\P$ is countably closed
to find a stronger condition $t'$ that decides $\dot A\intersect t$. Further, we can find a stronger $t''$ such that for every node $a\in t$ there is
a node $b\in t''$ comparable to it such that $t''\forces b\in\dot A$. Iterating this, in what we call the {\it bootstrap} argument, we build a
descending sequence $t_0> t_1>\cdots$ in $\P$ such that $t_{n+1}$ decides $\dot A\intersect t_n$ and $t_{n+1}$ forces that every element of $t_n$ is
comparable to a node in $\dot A\intersect t_{n+1}$. It follows that the limit tree $t_\omega=\bigcup t_n$ decides $\dot A\intersect t_\omega$ and
forces that it is a maximal antichain in $t_\omega$. For each node $a\in t$, let $b_a$ be a branch through $t_\omega$ containing $a$ and passing
through an element of the set $\dot A\intersect t_\omega$ decided by $t_\omega$. Let $\bar t$ be $t\union\{ b_a \mid a\in t_\omega\}$ be the
resulting tree. This is a condition in $\P$, stronger than $t_\omega$, but any node on the top level of $\bar t$ is comparable to something in $\dot
A\intersect t_\omega$. Thus, $\bar t$ forces that no additional nodes can be added to the antichain $\dot A$ above the height of $\bar t$: the
antichain is ``sealed''. Thus, $\bar t$ forces that $\dot A$ is countable, and so $T$ is Souslin.

Next, we show that $T$ has the unique branch property. We have to show that forcing with $T$ adds a unique branch through $T$. The combined forcing
that we are considering is $\P*\dot T$, where we first add the tree and then force with it. Given any condition $\<t,\dot q>$, we may strengthen $t$
to a tree of successor height deciding the particular value of $\dot q$, and then strengthen that value to a node on the maximal level of $t$. Thus,
this two-step forcing has a dense set $D$ of conditions of the form $\<t,\check q>$, where $t$ is a normal $(\alpha+1)$-tree for some countable
ordinal $\alpha$ and $q$ is a node on the $\alpha^{\rm th}$ level of $t$. The point is that forcing with $\P*\dot T$ is equivalent to forcing with
$D$, since $D$ is dense, but $D$ has the advantage of being countably closed. Suppose that $\tau$ is a $D$-name for another cofinal branch through
$T$, and this is forced by some condition $\<t_0,p_0>\in D$. Extend this to a stronger condition $\<t_1,p_1>\in D$ that decides $\tau\intersect t_0$,
and so on to build $\<t_n,p_n>>\<t_{n+1},p_{n+1}>$, with the next condition deciding $\tau$ on the previous tree. Let $t_\omega=\bigcup t_n$ be the
limit tree and $p_\omega=\bigcup p_n$ be the limit of the nodes $p_n$ in this tree. Let $b=\tau_{\union\<t_n,p_n>}$ be the branch through $t_\omega$
decided by the conditions $\<t_n,p_n>$. Extend $t_\omega$ to a tree $\bar t$ by adding branches through every node, to form a maximal level, but
without adding the branch $b$. Thus, $\<\bar t,p_\omega>$ is a condition in $D$, stronger than $\<t_0,p_0>$, forcing that $\tau$ cannot continue past
this level, contrary to our assumption that $\tau$ was forced to name a new cofinal branch.

The proof that $T$ has the absolute unique branch property is similar. Now, we are forcing with $\P*\dot T*\dot T$, and there is a dense set $D$ of
conditions of the form $\<t,p_1,p_2>$, where $t$ is a normal $(\alpha+1)$-tree for some countable $\alpha$ and $p_1$ and $p_2$ are two nodes in
$t(\alpha)$. If $\tau$ named a new cofinal branch, one could carry out the previous argument to build a condition deciding this branch, but then
extend the tree so as to prevent this new branch from extending. \qed

A similar argument, going beyond the notions considered in Figure \ref{fig:ImplicationDiagram}, shows that a generic Souslin tree $T$ is {\it
$n$-absolutely UBP} for every natural number $n$, meaning that forcing with $T^{n+1}$ adds exactly $n+1$ new cofinal branches through $T$. In fact,
generic Souslin trees exhibit even stronger rigidity properties:

\begin{definition}\rm
\label{definition:SouslinOffTheGenericBranch} A Souslin tree $T$ is {\it Souslin off the generic branch} if after forcing with $T$ to add a generic
branch $b$, then $T_p$ remains Souslin for every node $p$ not on $b$. More generally, $T$ has property $P$ {\it off the generic branch} if after
forcing with $T$ to add a generic branch $b$, then every $T_p$ has property $P$ for $p\notin b$.  We say that $T$ is {\it
  $n$-fold Souslin off the generic branch} if after forcing with
$T^n$, which adds $n$ branches $b_1,\ldots,b_n$, then $T_p$ remains Souslin for any $p$ not on any $b_i$.
\end{definition}%
Thus, another way to say that $T$ has the unique branch property is to say that it is Aronszajn off the generic branch. In particular, if a tree is
Souslin off the generic branch, then it will have the unique branch property, since a second branch would betray the Souslin-ness of that part of the
tree. Analogously, if a tree is $n$-fold Souslin off the generic branch, it must be $n$-absolutely UBP.

The other notions of rigidity that we introduced can be strengthened analogously. Given a property $P$, we say that $T$ is $n$-absolutely $P$ if $T$
has property $P$ in every generic extension obtained by forcing with $T^n$. Then the remarks about the preservation of $\omega_1$ carry over. Thus,
if a tree $T$ is $n$-absolutely rigid, then forcing with $T^n$ preserves $\omega_1$, and if $T$ is $n$-absolutely UBP, then forcing with $T^{n+1}$
preserves $\omega_1$.

\begin{thm}
Every $V$-generic Souslin tree is Souslin off the generic branch. Indeed, such trees are $n$-fold Souslin off the generic branch for every natural
number $n$. Consequently, they are also $n$-absolutely UBP for every $n<\omega$, and $n$-absolutely rigid, and so
on.\label{thm:GenericSouslinTreeIsSouslinOffGenericBranch}
\end{thm}

\proof Suppose that we force with $\P$ to add a $V$-generic Souslin tree $T$ and then force with the tree $T$ itself. This two-step forcing is
equivalent to the forcing $\P*\dot T$, where $\dot T$ is the $\P$-name for the resulting generic tree added by $\P$. We have seen already that this
forcing has a dense set $D$ consisting of normal $(\alpha+1)$-trees $t$ and nodes $b$ on the $\alpha^{\rm th}$ level $t(\alpha)$. This is countably
closed forcing. Suppose $\<t,q>$ is such a condition in $D$ and $p$ is a node in $t$ incomparable to $q$. We claim that the resulting tree $T_p$ is a
Souslin tree in the extension $V[T,b]$. To see this, suppose that $\dot A$ is a $D$-name for an antichain in $T_p$. The bootstrap argument of Theorem
\ref{thm:GenericSouslinTreesAreRigid} shows that any condition in $D$ can be extended to a stronger condition forcing that $\dot A$ is bounded.
Hence, $T_p$ is Souslin in $V[T,b]$, and so $T$ is Souslin off the generic branch.

The $n$-fold version is similar. Forcing with $\P$ to add a Souslin tree $T$ and then forcing with $T^n$ is the same as forcing with $\P*(\dot T)^n$.
This poset has a dense set $D$ that is isomorphic to the set of conditions of the form $\<t,q_1,\ldots,q_n>$, where $t$ is a normal $(\alpha+1)$-tree
for some countable ordinal $\alpha$ and $q_1,\ldots,q_n$ are $n$ many nodes on the $\alpha^{\rm th}$ level $t(\alpha)$. This forcing is countably
closed. If $p$ is a node in $t$ and incomparable to each $q_i$, then the argument of the previous paragraph shows that the resulting subtree $T_p$
will still be Souslin in $V[T,b_1,\ldots,b_n]$, as desired. \qed

Let us now turn to the construction of Souslin trees from $\diamondsuit$. This is the combinatorial principle asserting that there is a sequence
$\<D_\alpha\mid\alpha<\omega_1>$ such that for every $A\of\omega_1$ the set $\{\alpha\mid A\intersect\alpha=D_\alpha\}$ is stationary.

\begin{thm}[Jensen] \label{thm:DiamondImpliesRigidSuslinTree}
If\/ $\diamondsuit$ holds, then there is a rigid Souslin tree.
\end{thm}

\proof This construction is widely known (see \cite{Jech:SetTheory3rdEdition} or \cite{TSP}), but we include it as a warm-up to Theorem
\ref{thm:DiamondImpliesnFoldSouslin} and the arguments of Sections \ref{Section:SeparatingTheRigidityNotions} and \ref{Section:AutoTowers}. Fix any
$\diamondsuit$ sequence $\<D_\alpha\mid \alpha<\omega_1>$, using countable ordinals for the nodes. We first explain merely how to construct a Souslin
tree $T$, by recursively constructing the levels $T(\alpha)$. We begin with a sole root node. At successor stages, if the current top level
$T(\alpha)$ is defined, then we give each of these nodes two immediate successors in $T(\alpha+1)$. The nontrivial case occurs when $T|\alpha$ is
defined and $\alpha$ is a limit ordinal. If $D_\alpha$ happens to be a maximal antichain in $T|\alpha$, then for every node $p\in T|\alpha$, select a
branch $b_p\in[T |\alpha]$ with $p\in b_p$ and $b_p\cap D_\alpha\neq\emptyset$. For each such branch $b_p$, place a node in $T(\alpha)$ on top of it.
In this case, we say that $T(\alpha)$ has {\it sealed} the antichain $D_\alpha$. If $D_\alpha$ is not a maximal antichain in $T|\alpha$, then carry
out the same construction, but dropping the requirement involving $D_\alpha$. The resulting tree $T$ is clearly a normal $\omega_1$-tree. To see that
it is Souslin, suppose that $A\subset T$ is a maximal antichain. By some simple closure arguments, there is a club of $\alpha$ such that
$A\cap\alpha$ is a maximal antichain in $T|\alpha$. By $\diamondsuit$, the set of $\alpha$ with $A\cap\alpha=D_\alpha$ is stationary, and so there is
an $\alpha<\omega_1$ such that $D_\alpha=A\cap\alpha$ is a maximal antichain in $T|\alpha$. In this case, we sealed the antichain, and every node in
$T(\alpha)$ lies above an element of $A\intersect (T|\alpha)$. Therefore, all elements of $T$ above level $\alpha$ are comparable to an element of
$A\intersect (T|\alpha)$, and so $A=A\intersect (T|\alpha)$ is countable, as desired.

To ensure that $T$ is rigid, we now fold an additional step into the construction. At a limit stage $\alpha$, if $D_\alpha$ happens to code a
nontrivial automorphism $f$ of the tree $T|\alpha$, then we find a branch $b\in[T|\alpha]$ such that $f[b]\neq b$. For each $p\in T|\alpha$, select a
branch $b_p\in[T_p|\alpha]$ with $b_p\neq f[b]$, and let $T(\alpha)=\{ b_p\mid p\in T|\alpha\}\union\{b\}$. Note that because we have included $b$
but not $f[b]$ in $T(\alpha)$, we have prevented $f$ from being the initial segment of an automorphism of $T|(\alpha+1)$. It follows that the final
tree $T$ is rigid, because if $f:T\cong T$ were a nontrivial automorphism of the full tree, then there would be a stationary set of $\alpha$ such
that $D_\alpha$ codes $f\restrict T|\alpha$, and at such a stage, once we are into the nontrivial part of $f$, we specifically constructed the
$\alpha^{\rm th}$ level so as to exclude the possibility that $f\restrict (T|\alpha)$ could be extended to an automorphism of $T|(\alpha+1)$. Since
this contradicts our assumption that $f$ was an automorphism of $T$, we have constructed a rigid Souslin tree. \qed

It was also shown by Jensen that $\diamondsuit$ implies the existence of Souslin trees that fail to be rigid, in a very strong way, by satisfying
strong forms of homogeneity (see \cite{TSP}).

The key method in the proof of Theorem \ref{thm:DiamondImpliesRigidSuslinTree} was to anticipate via $\diamondsuit$ the potential antichains and
automorphisms of the tree, and seal them into a level of the tree, preventing them from growing into full counterexamples. An elaboration of this
method allows us also to attain the higher degrees of rigidity.

\begin{thm}\label{thm:DiamondImpliesnFoldSouslin}
If\/ $\diamondsuit$ holds, then there is a Souslin tree that is $n$-fold Souslin off the generic branch, for every finite $n$, and consequently also
$n$-absolutely UBP and $n$-absolutely rigid and so on.
\end{thm}

\goodbreak
\proof We again construct the tree $T$ by induction on the levels $T|\alpha$, using a fixed $\diamondsuit$ sequence $\<D_\alpha\mid\alpha<\omega_1>$.
We begin as before with a sole root node, and at successors we give every node two immediate successors. So assume that $\alpha$ is a limit ordinal
and let $t=T|\alpha$ be the tree constructed up to stage $\alpha$. While we could anticipate and then seal potential antichains and potential
automorphisms, it will be unnecessary to do so explicitly. Rather, it will suffice to anticipate and then seal the potential violations to the tree
being $n$-fold Souslin off the generic branch. Specifically, suppose that for some natural number $n$ and some $\vec q\in t^n$, the set $D_\alpha$
codes a $t^n_{[\vec q]}$-name $\dot A$ for a maximal antichain in $t_p$ for some $p\perp q_i$. More precisely, associated to this name is the
function $\vec r\mapsto \dot A_{\vec r}$, where $\dot A_{\vec r}$ is the set of elements in $t_p$ that the condition $\vec r$ determines to be in the
name $\dot A$. We assume that every $\dot A_{\vec r}$ is an antichain in $t_p$; if $\vec s$ extends $\vec r$, then $\dot A_{\vec r}\of \dot A_{\vec
s}$; and for any $q\in t_p$ and any $\vec r$ extending $\vec q$, there is some $\vec s$ extending $\vec r$ such that $\dot A_{\vec s}$ contains an
element compatible with $q$ in $t_p$. Using these properties and the fact that $t_p$ is countable, we may successively extend $\vec q$ to meet a
certain countable collection of dense sets so as to build a cofinal branch $\vec b\in [t^n]$ extending $\vec q$ so that the corresponding combined
information $\dot A_{\vec b}=\Union\{\dot A_{\vec b\restrict\beta}\mid\beta<\alpha\}$ determined by the nodes on this branch is a maximal antichain
in $t_p$. Since $p\perp q_i$, the node $p$ is not on any of the branches $b_i$. Now, for every node $q$ in $t_{[p]}$, we find a branch
$b_q\in[t_{[q]}]$ going through the maximal antichain $\dot A_{\vec b}$. For each node $q\in t$ incomparable with $p$, we find any branch
$b_q\in[t_{[q]}]$. Let $B$ be the resulting combined set of branches, including $\vec b$ and all $b_q$ for $q\in t$, and define $T(\alpha)$ to have
nodes above exactly the branches in $B$. By construction, this is a normal $(\alpha+1)$-tree.

We now argue for every finite $n$ that the resulting tree $T$ is $n$-fold Souslin off the generic branch. Suppose for some $n$, some condition $\vec
q=\<q_1,\ldots,q_n>\in T^n$ and some $p\perp q_i$, there is a $T^n$-name $\dot A$ such that $\vec q$ forces that $\dot A$ is a maximal antichain in
$T_p$. The name $\dot A$ is coded by a subset $\dot A_0\of\omega_1$. There is a club of $\alpha$ such that $T\intersect\alpha=T|\alpha$, and $\dot
A_0\intersect\alpha$ is the code for $\dot A\intersect(T|\alpha)$, and such that this is actually a $(T|\alpha)^n$-name for a maximal antichain in
$(T|\alpha)_p$. By $\diamondsuit$, therefore, there is some $\alpha$ like this such that in addition, $D_\alpha$ is the code of $\dot A\intersect
(T|\alpha)$. Therefore, our construction above exactly ensures that $T(\alpha)$ has nodes $\vec b=\<b_1,\ldots,b_n>$ that already determine $\dot
A\intersect (T_p|\alpha)$ to be a maximal antichain in $T_p|\alpha$, and furthermore, that every node in $T_p(\alpha)$ passed through a node in the
set $\dot A_{\vec b}=\dot A\intersect (T_p|\alpha)$ determined by $\vec b$. Thus, $\vec b$ forces via $T^n$ that $\dot A$ has been sealed by level
$\alpha$, and consequently that $\dot A=\dot A_{\vec b}$. So $\vec q$ did not force that $\dot A$ was unbounded. Consequently, since $n$, $p$ and
$\vec q$ were arbitrary, we have established that $T$ is $n$-fold Souslin off the generic branch for every $n$. It follows that $T$ is Souslin, and
rigid, and $n$-absolutely rigid, and $n$-absolutely UBP and so on. \qed

We close this section by remarking that we have checked in detail that Jensen's original tree (see the construction of \cite[Thm. V.1.1]{TSP}) is
also $n$-fold Souslin off the generic branch.

\section{Separating the rigidity notions}
\label{Section:SeparatingTheRigidityNotions}

We turn now to our main result, proving that the diagram of implications in Figure \ref{fig:ImplicationDiagram} is complete and omits no
\ZFC-provable implication. For convenience, we reproduce the diagram here.
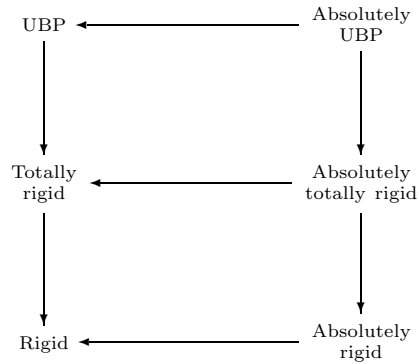
\begin{figure}[htb]
\begin{diagram}[width=6em,height=3em,objectstyle=\rm,textflow]
 {\scriptstyle UBP}             &   \lTo        &  {Absolutely\atop UBP}               \\
 \dTo            &               &  \dTo                         \\
 {Totally\atop rigid}   &   \lTo        &  {Absolutely\atop totally\ rigid}     \\
  \dTo           &               &  \dTo                         \\
 {\scriptstyle Rigid}           &  \lTo         &  {Absolutely\atop rigid}             \\
\end{diagram}
\end{figure}
\begin{obs}\label{obs:SufficesToFind}
In order to show that the implication diagram above is complete, it suffices to find
\begin{enumerate}
 \item An absolutely rigid tree that is not totally rigid,
 \item An absolutely totally rigid tree that is not UBP, and
 \item A UBP tree that is not absolutely rigid.
\end{enumerate}
\end{obs}

\proof Assume there are trees exhibiting 1, 2 and 3.  A tree as in 1 shows that there are no implications from the bottom row to the middle row. A
tree as in 2 shows that there are no implications from the middle row to the top row. And a tree as in 3 shows that there are no implications from
the left column to the right column. This refutes any potential implication either going (in any sense) up or to the right. It follows that there are
no missing implications, because the transitive closure of the diagram already exhibits all implications going down and to the left. \qed

So we will prove Theorem \ref{thm:MainTheorem} by finding trees fulfilling the requirements of Observation \ref{obs:SufficesToFind}.

\begin{thm}\label{thm:AbsolutelyRigidNotTotallyRigid}
If there is an absolutely rigid Souslin tree, then there is an absolutely rigid Souslin tree that is not totally rigid.
\end{thm}

This is a consequence of the following more general lemma.

\goodbreak
\begin{lem}
\label{separation} If $T$ is a tree, then there is a tree $C(T)$ such that:
\begin{enumerate}
  \item $T$ is rigid if and only if $C(T)$ is rigid.
  \item $T$ is absolutely rigid if and only if $C(T)$ is absolutely rigid.
  \item $T$ is Souslin if and only if $C(T)$ is Souslin.
  \item $C(T)$ is not totally rigid.
  \item $T$ and $C(T)$ are forcing equivalent.
\end{enumerate}
\end{lem}

\begin{figure}[htb]
\label{fig:gluedTree}
\begin{center}
\setlength{\unitlength}{0.00043333in}
\begingroup\makeatletter\ifx\SetFigFont\undefined%
\gdef\SetFigFont#1#2#3#4#5{%
  \reset@font\fontsize{#1}{#2pt}%
  \fontfamily{#3}\fontseries{#4}\fontshape{#5}%
  \selectfont}%
\fi\endgroup%
{\renewcommand{\dashlinestretch}{30}
\begin{picture}(7295,4804)(0,-10)
\path(2412,3420)(3012,2220)(3612,3420)
\path(2412,3420)(2412,3421)(2413,3423)
	(2415,3427)(2418,3433)(2422,3441)
	(2427,3450)(2433,3461)(2440,3474)
	(2449,3487)(2459,3500)(2470,3513)
	(2483,3525)(2497,3537)(2513,3548)
	(2532,3557)(2553,3565)(2577,3571)
	(2605,3575)(2637,3577)(2673,3575)
	(2712,3570)(2750,3562)(2788,3552)
	(2824,3541)(2857,3528)(2887,3514)
	(2913,3501)(2935,3487)(2955,3473)
	(2971,3460)(2986,3446)(2999,3433)
	(3012,3420)(3025,3407)(3038,3394)
	(3053,3380)(3069,3367)(3089,3353)
	(3111,3339)(3137,3326)(3167,3312)
	(3200,3299)(3236,3288)(3274,3278)
	(3312,3270)(3351,3265)(3387,3263)
	(3419,3265)(3447,3269)(3471,3275)
	(3492,3283)(3511,3292)(3527,3303)
	(3541,3315)(3554,3327)(3565,3340)
	(3575,3353)(3584,3366)(3591,3379)
	(3597,3390)(3602,3399)(3606,3407)
	(3609,3413)(3611,3417)(3612,3419)(3612,3420)
\path(3012,2220)(3612,1620)(4812,2220)
\path(1212,2820)(1812,1620)(2412,2820)
\path(1212,2820)(1212,2821)(1213,2823)
	(1215,2827)(1218,2833)(1222,2841)
	(1227,2850)(1233,2861)(1240,2874)
	(1249,2887)(1259,2900)(1270,2913)
	(1283,2925)(1297,2937)(1313,2948)
	(1332,2957)(1353,2965)(1377,2971)
	(1405,2975)(1437,2977)(1473,2975)
	(1512,2970)(1550,2962)(1588,2952)
	(1624,2941)(1657,2928)(1687,2914)
	(1713,2901)(1735,2887)(1755,2873)
	(1771,2860)(1786,2846)(1799,2833)
	(1812,2820)(1825,2807)(1838,2794)
	(1853,2780)(1869,2767)(1889,2753)
	(1911,2739)(1937,2726)(1967,2712)
	(2000,2699)(2036,2688)(2074,2678)
	(2112,2670)(2151,2665)(2187,2663)
	(2219,2665)(2247,2669)(2271,2675)
	(2292,2683)(2311,2692)(2327,2703)
	(2341,2715)(2354,2727)(2365,2740)
	(2375,2753)(2384,2766)(2391,2779)
	(2397,2790)(2402,2799)(2406,2807)
	(2409,2813)(2411,2817)(2412,2819)(2412,2820)
\path(1812,1620)(2412,1020)(3612,1620)
\path(3612,4020)(4212,2820)(4812,4020)
\path(3612,4020)(3612,4021)(3613,4023)
	(3615,4027)(3618,4033)(3622,4041)
	(3627,4050)(3633,4061)(3640,4074)
	(3649,4087)(3659,4100)(3670,4113)
	(3683,4125)(3697,4137)(3713,4148)
	(3732,4157)(3753,4165)(3777,4171)
	(3805,4175)(3837,4177)(3873,4175)
	(3912,4170)(3950,4162)(3988,4152)
	(4024,4141)(4057,4128)(4087,4114)
	(4113,4101)(4135,4087)(4155,4073)
	(4171,4060)(4186,4046)(4199,4033)
	(4212,4020)(4225,4007)(4238,3994)
	(4253,3980)(4269,3967)(4289,3953)
	(4311,3939)(4337,3926)(4367,3912)
	(4400,3899)(4436,3888)(4474,3878)
	(4512,3870)(4551,3865)(4587,3863)
	(4619,3865)(4647,3869)(4671,3875)
	(4692,3883)(4711,3892)(4727,3903)
	(4741,3915)(4754,3927)(4765,3940)
	(4775,3953)(4784,3966)(4791,3979)
	(4797,3990)(4802,3999)(4806,4007)
	(4809,4013)(4811,4017)(4812,4019)(4812,4020)
\path(4212,2820)(4812,2220)(6012,2820)
\put(7212,3420){\ellipse{150}{150}}
\put(2412,1020){\blacken\ellipse{150}{150}}
\put(2412,1020){\ellipse{150}{150}}
\put(3612,1620){\blacken\ellipse{150}{150}}
\put(3612,1620){\ellipse{150}{150}}
\put(4812,2220){\blacken\ellipse{150}{150}}
\put(4812,2220){\ellipse{150}{150}}
\put(6012,2820){\blacken\ellipse{150}{150}}
\put(6012,2820){\ellipse{150}{150}}
\put(1212,420){\blacken\ellipse{150}{150}}
\put(1212,420){\ellipse{150}{150}}
\path(6012,2820)(5412,3420)
\path(4812,4620)(5412,3420)(6012,4620)
\path(12,2220)(612,1020)(1212,2220)
\path(612,1020)(1212,420)(2412,1020)
\dashline{45.000}(6012,2820)(7212,3420)
\path(12,2220)(12,2221)(13,2223)
	(15,2227)(18,2233)(22,2241)
	(27,2250)(33,2261)(40,2274)
	(49,2287)(59,2300)(70,2313)
	(83,2325)(97,2337)(113,2348)
	(132,2357)(153,2365)(177,2371)
	(205,2375)(237,2377)(273,2375)
	(312,2370)(350,2362)(388,2352)
	(424,2341)(457,2328)(487,2314)
	(513,2301)(535,2287)(555,2273)
	(571,2260)(586,2246)(599,2233)
	(612,2220)(625,2207)(638,2194)
	(653,2180)(669,2167)(689,2153)
	(711,2139)(737,2126)(767,2112)
	(800,2099)(836,2088)(874,2078)
	(912,2070)(951,2065)(987,2063)
	(1019,2065)(1047,2069)(1071,2075)
	(1092,2083)(1111,2092)(1127,2103)
	(1141,2115)(1154,2127)(1165,2140)
	(1175,2153)(1184,2166)(1191,2179)
	(1197,2190)(1202,2199)(1206,2207)
	(1209,2213)(1211,2217)(1212,2219)(1212,2220)
\path(4812,4620)(4812,4621)(4813,4623)
	(4815,4627)(4818,4633)(4822,4641)
	(4827,4650)(4833,4661)(4840,4674)
	(4849,4687)(4859,4700)(4870,4713)
	(4883,4725)(4897,4737)(4913,4748)
	(4932,4757)(4953,4765)(4977,4771)
	(5005,4775)(5037,4777)(5073,4775)
	(5112,4770)(5150,4762)(5188,4752)
	(5224,4741)(5257,4728)(5287,4714)
	(5313,4701)(5335,4687)(5355,4673)
	(5371,4660)(5386,4646)(5399,4633)
	(5412,4620)(5425,4607)(5438,4594)
	(5453,4580)(5469,4567)(5489,4553)
	(5511,4539)(5537,4526)(5567,4512)
	(5600,4499)(5636,4488)(5674,4478)
	(5712,4470)(5751,4465)(5787,4463)
	(5819,4465)(5847,4469)(5871,4475)
	(5892,4483)(5911,4492)(5927,4503)
	(5941,4515)(5954,4527)(5965,4540)
	(5975,4553)(5984,4566)(5991,4579)
	(5997,4590)(6002,4599)(6006,4607)
	(6009,4613)(6011,4617)(6012,4619)(6012,4620)
\put(612,1620){\makebox(0,0)[b]{{\SetFigFont{12}{14.4}{\familydefault}{\mddefault}{\updefault}$T_0$}}}
\put(1812,2220){\makebox(0,0)[b]{{\SetFigFont{12}{14.4}{\familydefault}{\mddefault}{\updefault}$T_1$}}}
\put(3012,2820){\makebox(0,0)[b]{{\SetFigFont{12}{14.4}{\familydefault}{\mddefault}{\updefault}$T_2$}}}
\put(4212,3420){\makebox(0,0)[b]{{\SetFigFont{12}{14.4}{\familydefault}{\mddefault}{\updefault}$T_3$}}}
\put(5412,4020){\makebox(0,0)[b]{{\SetFigFont{12}{14.4}{\familydefault}{\mddefault}{\updefault}$T_4$}}}
\put(1212,45){\makebox(0,0)[b]{{\SetFigFont{12}{14.4}{\rmdefault}{\mddefault}{\updefault}$r_0$}}}
\put(2412,645){\makebox(0,0)[b]{{\SetFigFont{12}{14.4}{\rmdefault}{\mddefault}{\updefault}$r_1$}}}
\put(3612,1245){\makebox(0,0)[b]{{\SetFigFont{12}{14.4}{\rmdefault}{\mddefault}{\updefault}$r_2$}}}
\put(4812,1845){\makebox(0,0)[b]{{\SetFigFont{12}{14.4}{\rmdefault}{\mddefault}{\updefault}$r_3$}}}
\put(6012,2445){\makebox(0,0)[b]{{\SetFigFont{12}{14.4}{\rmdefault}{\mddefault}{\updefault}$r_4$}}}
\end{picture}
}
\end{center}
\caption{The tree $C(T)$ has $\omega$ many copies of $T$} \label{fig:C-of-T}
\end{figure}
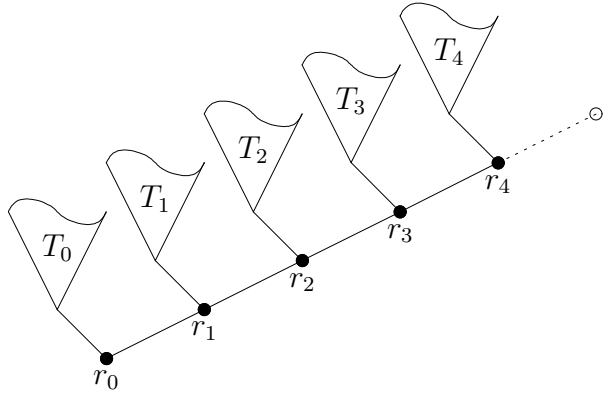

\proof The tree $C(T)$ is built by gluing $\omega$ copies of $T$ together as in figure \ref{fig:C-of-T}.
Denoting the $n^{\rm th}$ copy of $T$ by $T_n$, we add the new nodes $r_1,r_2,\ldots$, having no limit in the tree, and use them to glue the copies
of $T$ together in such a way that $r_0$ is the root of the new tree, and the immediate successors of $r_n$ are the root of $T_n$ and $r_{n+1}$.
Clearly, antichains in $C(T)$ consist essentially of copies of $\omega$ many antichains in $T$, and so 3 holds by the pigeonhole principle. The tree
$C(T)$ is not totally rigid, since it is built from copies of $T$, and so 4 holds. Forcing with $C(T)$ amounts to selecting a copy of $T$ and forcing
with it, and so 5 holds.

Consider statement 1. Clearly any nontrivial automorphism of $T$ gives rise to a nontrivial automorphism of $C(T)$, so $T$ is rigid if $C(T)$ is.
Conversely, suppose $T$ is rigid and $\pi$ is a non-trivial automorphism of $C(T)$. Since each $T_n$ is a copy of $T$ and hence rigid, it follows
easily that $\pi$ must move the root of some $T_n$ to the adjacent $r_{n+1}$. Thus, $\pi$ witnesses that $T_n$ is isomorphic to $C(T)_{r_{n+1}}$. By
construction, however, $C(T)_{r_{n+1}}$ is isomorphic to $C(T)$. We conclude that $T_n$, and hence $T$, is isomorphic to $C(T)$. But this contradicts
our assumptions that $T$ was rigid, while $C(T)$ was not. So 1 holds.

Lastly, consider statement 2. By 5, forcing with $T$ or $C(T)$ gives rise to the same extension. Applying 1 in such an extension, since the
definition of $C(T)$ is absolute, we conclude that $T$ is rigid there if and only if $C(T)$ is rigid. Thus, in the ground model, $T$ was absolutely
rigid if and only if $C(T)$ was absolutely rigid. So 2 holds. \qed

This proves Theorem \ref{thm:AbsolutelyRigidNotTotallyRigid} and fulfills the first requirement of Observation \ref{obs:SufficesToFind}. We turn now
to the second requirement.

\begin{thm}\label{thm:AbsolutelyTotallyRigidNotUBP}
If there is a (Souslin) tree with the absolute unique branch property,
then there is an absolutely totally rigid (Souslin) tree without the unique branch property.
%
\end{thm}

To prove this theorem, we introduce another notion of (non)rigidity: a structure is \emph{Hopfian} if it is isomorphic to a proper substructure of
itself. The following Lemma answers a question by Martin Weese (private communication).

\begin{lem}
\label{lem:UBP-implies-not-Hopfian} Normal UBP trees are not Hopfian.
\end{lem}

\proof Assume the contrary. Suppose $T$ is a normal UBP tree, $A\rsubset  T$ and $\pi:T\cong T\rest A$. We claim that there is a $p\in T$ such that
$\pi(p)\perp p$. Such a $p$ forces via $T$ that there are at least two new cofinal branches, for if $b$ is a generic branch of $T$ through $p$, then
$\pi[b]$ determines a cofinal branch through $\pi(p)$, and as $p\perp\pi(p)$, these branches are different. This contradicts the unique branch
property. To find such a $p$, let $\bp\in  T\ohne A$. If $\pi(\bp)\perp\bp$, we're done. Otherwise, $\pi(\bp)>\bp$. By the normality of $T$, let
$p\in T$ be such that $p>\bp$ and $p\perp\pi(\bp)$. Then $\pi(p)>\pi(\bp)$, and so $p\perp\pi(p)$, because $p\perp\pi(\bp)$, as desired. \qed

\begin{lem}\label{lem:totally-rigid-not-UBP} If\/ $T$ is a normal $\omega_1$-tree, then there is a normal
$\omega_1$-tree $C'(T)$ such that:
\begin{enumerate}
  \item \label{item:not-UBP} $C'(T)$ is not UBP.
  \item \label{item:totally-rigid} If\/ $T$ is UBP, then $C'(T)$
    is totally rigid.
  \item \label{item:abs-totally-rigid} If\/ $T$ is absolutely
    UBP, then $C'(T)$ is absolutely totally rigid.
  \item \label{item:souslinness-preserved} $T$ is Souslin if and only if $C'(T)$ is Souslin.
  \item \label{item:omega-splitting} $C'(T)$ is $\omega$-splitting.
\end{enumerate}
\end{lem}

\begin{figure}[htb]\label{fig:TotallyRigidNotUBPTree}
\setlength{\unitlength}{0.0004in}
\begingroup\makeatletter\ifx\SetFigFont\undefined%
\gdef\SetFigFont#1#2#3#4#5{%
  \reset@font\fontsize{#1}{#2pt}%
  \fontfamily{#3}\fontseries{#4}\fontshape{#5}%
  \selectfont}%
\fi\endgroup%
{\renewcommand{\dashlinestretch}{30}
\begin{picture}(11424,3489)(0,-10)
\put(11412,1812){\makebox(0,0)[b]{\smash{{{\SetFigFont{12}{14.4}{\rmdefault}{\mddefault}{\updefault}$\cdots$}}}}}
\put(8637,2412){\makebox(0,0)[b]{\smash{{{\SetFigFont{12}{14.4}{\rmdefault}{\mddefault}{\updefault}$T|C_3$}}}}}
\put(3462,2412){\makebox(0,0)[b]{\smash{{{\SetFigFont{12}{14.4}{\rmdefault}{\mddefault}{\updefault}$T|C_1$}}}}}
\put(837,2412){\makebox(0,0)[b]{\smash{{{\SetFigFont{12}{14.4}{\rmdefault}{\mddefault}{\updefault}$T|C_0$}}}}}
\put(6087,2412){\makebox(0,0)[b]{\smash{{{\SetFigFont{12}{14.4}{\rmdefault}{\mddefault}{\updefault}$T|C_2$}}}}}
\path(2604,3354)(2606,3355)(2609,3356)
	(2615,3359)(2624,3364)(2636,3369)
	(2652,3376)(2670,3384)(2690,3392)
	(2712,3401)(2736,3410)(2761,3419)
	(2788,3427)(2816,3435)(2846,3442)
	(2879,3448)(2914,3454)(2952,3458)
	(2993,3461)(3036,3462)(3079,3461)
	(3118,3459)(3152,3455)(3181,3451)
	(3203,3447)(3220,3443)(3232,3438)
	(3241,3434)(3247,3430)(3252,3426)
	(3257,3422)(3263,3417)(3272,3412)
	(3284,3406)(3301,3400)(3323,3393)
	(3352,3384)(3386,3375)(3425,3365)
	(3468,3354)(3511,3343)(3550,3333)
	(3584,3324)(3613,3315)(3635,3308)
	(3652,3302)(3664,3296)(3673,3291)
	(3679,3286)(3684,3282)(3689,3278)
	(3695,3274)(3704,3270)(3716,3265)
	(3733,3261)(3755,3257)(3784,3253)
	(3818,3249)(3857,3247)(3900,3246)
	(3943,3247)(3984,3250)(4022,3254)
	(4057,3260)(4090,3266)(4120,3273)
	(4148,3281)(4175,3289)(4200,3298)
	(4224,3307)(4246,3316)(4266,3324)
	(4284,3332)(4300,3339)(4312,3344)
	(4321,3349)(4327,3352)(4330,3353)(4332,3354)
\path(12,3354)(14,3355)(17,3356)
	(23,3359)(32,3364)(44,3369)
	(60,3376)(78,3384)(98,3392)
	(120,3401)(144,3410)(169,3419)
	(196,3427)(224,3435)(254,3442)
	(287,3448)(322,3454)(360,3458)
	(401,3461)(444,3462)(487,3461)
	(526,3459)(560,3455)(589,3451)
	(611,3447)(628,3443)(640,3438)
	(649,3434)(655,3430)(660,3426)
	(665,3422)(671,3417)(680,3412)
	(692,3406)(709,3400)(731,3393)
	(760,3384)(794,3375)(833,3365)
	(876,3354)(919,3343)(958,3333)
	(992,3324)(1021,3315)(1043,3308)
	(1060,3302)(1072,3296)(1081,3291)
	(1087,3286)(1092,3282)(1097,3278)
	(1103,3274)(1112,3270)(1124,3265)
	(1141,3261)(1163,3257)(1192,3253)
	(1226,3249)(1265,3247)(1308,3246)
	(1351,3247)(1392,3250)(1430,3254)
	(1465,3260)(1498,3266)(1528,3273)
	(1556,3281)(1583,3289)(1608,3298)
	(1632,3307)(1654,3316)(1674,3324)
	(1692,3332)(1708,3339)(1720,3344)
	(1729,3349)(1735,3352)(1738,3353)(1740,3354)
\path(7788,3354)(7790,3355)(7793,3356)
	(7799,3359)(7808,3364)(7820,3369)
	(7836,3376)(7854,3384)(7874,3392)
	(7896,3401)(7920,3410)(7945,3419)
	(7972,3427)(8000,3435)(8030,3442)
	(8063,3448)(8098,3454)(8136,3458)
	(8177,3461)(8220,3462)(8263,3461)
	(8302,3459)(8336,3455)(8365,3451)
	(8387,3447)(8404,3443)(8416,3438)
	(8425,3434)(8431,3430)(8436,3426)
	(8441,3422)(8447,3417)(8456,3412)
	(8468,3406)(8485,3400)(8507,3393)
	(8536,3384)(8570,3375)(8609,3365)
	(8652,3354)(8695,3343)(8734,3333)
	(8768,3324)(8797,3315)(8819,3308)
	(8836,3302)(8848,3296)(8857,3291)
	(8863,3286)(8868,3282)(8873,3278)
	(8879,3274)(8888,3270)(8900,3265)
	(8917,3261)(8939,3257)(8968,3253)
	(9002,3249)(9041,3247)(9084,3246)
	(9127,3247)(9168,3250)(9206,3254)
	(9241,3260)(9274,3266)(9304,3273)
	(9332,3281)(9359,3289)(9384,3298)
	(9408,3307)(9430,3316)(9450,3324)
	(9468,3332)(9484,3339)(9496,3344)
	(9505,3349)(9511,3352)(9514,3353)(9516,3354)
\path(5196,3354)(5198,3355)(5201,3356)
	(5207,3359)(5216,3364)(5228,3369)
	(5244,3376)(5262,3384)(5282,3392)
	(5304,3401)(5328,3410)(5353,3419)
	(5380,3427)(5408,3435)(5438,3442)
	(5471,3448)(5506,3454)(5544,3458)
	(5585,3461)(5628,3462)(5671,3461)
	(5710,3459)(5744,3455)(5773,3451)
	(5795,3447)(5812,3443)(5824,3438)
	(5833,3434)(5839,3430)(5844,3426)
	(5849,3422)(5855,3417)(5864,3412)
	(5876,3406)(5893,3400)(5915,3393)
	(5944,3384)(5978,3375)(6017,3365)
	(6060,3354)(6103,3343)(6142,3333)
	(6176,3324)(6205,3315)(6227,3308)
	(6244,3302)(6256,3296)(6265,3291)
	(6271,3286)(6276,3282)(6281,3278)
	(6287,3274)(6296,3270)(6308,3265)
	(6325,3261)(6347,3257)(6376,3253)
	(6410,3249)(6449,3247)(6492,3246)
	(6535,3247)(6576,3250)(6614,3254)
	(6649,3260)(6682,3266)(6712,3273)
	(6740,3281)(6767,3289)(6792,3298)
	(6816,3307)(6838,3316)(6858,3324)
	(6876,3332)(6892,3339)(6904,3344)
	(6913,3349)(6919,3352)(6922,3353)(6924,3354)
\path(5196,3354)(6060,762)(6924,3354)
\path(7212,12)(912,762)
\path(3462,762)(7212,12)
\dottedline{45}(7212,12)(11412,762)
\path(8637,762)(7212,12)
\path(6087,762)(7212,12)
\path(2604,3354)(3468,762)(4332,3354)
\path(12,3354)(876,762)(1740,3354)
\path(7788,3354)(8652,762)(9516,3354)
\end{picture}
}
\caption{The tree $C'(T)$}
\end{figure}
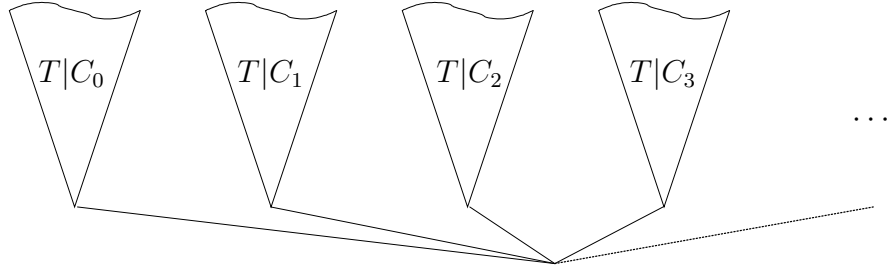

\proof Let $T$ be a normal $\omega_1$-tree. We define a descending sequence $\kla{C_n\st n<\omega}$ of club subsets of $\omega_1$ as follows. For any
subset $A$ of $\omega_1$, let $A'$ be the set of limit points of $A$ below $\omega_1$. Let $C_0=\{0\}\union\omega_1'$, and recursively define
$C_{n+1}=\{0\}\cup C_n'$. Let $C'(T)$ be the tree sketched in figure \ref{fig:TotallyRigidNotUBPTree}.
More precisely, $C'(T)$ consists of nodes $\{0\}\cup\bigcup_{n\in\omega}(\{n\}\times(T|C_n))$. Note that $T|C_n$ is an $\omega$-splitting normal
$\omega_1$ tree, since successive elements of $C_n$ jump over gaps of length at least $\omega$, and the limit nodes are unique since $C_n$ is closed.
It follows that $C'(T)$ is also a normal $\omega$-splitting tree, so statement \ref{item:omega-splitting} holds. Since uncountable antichains in $T$
give rise to uncountable antichains in $T|C$ and vice versa, statement \ref{item:souslinness-preserved} holds. Since every $T|C_n$ is dense in $T$,
it follows that forcing with $T$ is equivalent to forcing with $C'(T)$. Clearly, adding a branch to $T$ will add $\omega$ many branches to $C'(T)$,
so statement \ref{item:not-UBP} holds. Statement \ref{item:abs-totally-rigid} follows from statement \ref{item:totally-rigid}, since forcing with $T$
or $C'(T)$ gives rise to the same extensions, and so if $T$ has the unique branch property in such an extension, then by statement
\ref{item:totally-rigid} we know $C'(T)$ is also totally rigid there, as the definition of $C'(T)$ is absolute.

Lastly, consider statement \ref{item:totally-rigid}, and suppose that $T$ has the unique branch property. It follows easily that every $T|C_n$ also
has the unique branch property.\footnote{In fact, if $C$ is a club subset of $\omega_1$ containing
  $0$, then $T$ has the unique branch property if and only if $T|C$ also has it. So, the
  unique branch property is a very natural notion: In
  \cite{ForcingWithTrees}, it was shown that when investigating
  forcing extensions obtained by forcing with a Souslin tree $T$, the
  rigidity properties of $T|C$ matter, not those of $T$. I.e.,
  automorphisms of the complete Boolean algebra associated to $T$
  correspond to automorphisms of $T|C$, for some $C$ as above -- see
  \cite[Lemma 3.1.]{ForcingWithTrees}.}
By Lemma \ref{implications}, therefore, every $T|C_n$ is also totally rigid. Suppose towards contradiction that $C'(T)$ is not totally rigid. So
there is an isomorphism $\pi:C'(T)_\kla{i,\bp}\cong C'(T)_\kla{j,\bq}$ witnessing this. Since $T|C_i$ is totally rigid, it must be that $i\neq j$,
and we may assume $i<j$. We claim that there is an extension $\kla{i,p}$ of $\kla{i,\bp}$ whose image $\kla{j,q}=\pi(\kla{i,p})$ has $|p|_T=|q|_T$.
To see this, let $\vgamma$ and $\vgamma'$ be the continuous monotone enumerations of the clubs $C_i\ohne|\bp|_T$ and $C_j\ohne|\bq|_T$, respectively.
The fixed points of these enumerations form a club, and whenever $\alpha$ is such a common fixed point and $|r|_T=\alpha$, then
$|\pi(\kla{i,r})|_T=\alpha$ as well. So any $p>_T\bp$ at such a level $\alpha$ will have $|p|_T=|q|_T$, as we claimed. Since
$\pi(\kla{i,p})=\kla{j,q}$, the isomorphism $\pi$ induces an isomorphism $\pi':(T|C_i)_p\cong(T|C_j)_q$ on the underlying trees, defined by
$\pi'(r)=r'$ if and only if $\pi(\kla{i,r})=\kla{j,r'}$. If $p\neq q$, then $p\perp\pi'(p)$ and so forcing with the condition $p$ in $T|C_n$ will add
at least two new branches: the generic branch and the branch containing its image under $\pi'$. This contradicts the fact that $T|C_n$ is UBP.
Otherwise, we assume $p=q$, and so $\pi'$ is an isomorphism of $(T|C_i)_p$ with a proper subtree of itself. By Lemma
\ref{lem:UBP-implies-not-Hopfian}, this also contradicts the unique branch property of $(T|C_i)_p$. So statement \ref{item:totally-rigid} holds, and
the proof is complete. \qed

This proves Theorem \ref{thm:AbsolutelyTotallyRigidNotUBP} and fulfills the second requirement of Observation \ref{obs:SufficesToFind}. We turn now
to the third and most difficult requirement. For this, we will assume $\diamondsuit$ and construct a Souslin tree with the unique branch property,
but which is not absolutely rigid. The basic idea will be to construct a tree $T$ such that cofinal branches through $T$ automatically code
automorphisms of $T$. Thus, forcing with $T$ will necessarily add automorphisms to $T$. The difficulty will be to do this while retaining the unique
branch property. To assist with our construction, we introduce the concept of level-transitive group actions on trees.

\begin{definition}\rm
\label{def:level-transitive-actions} Let $T$ be a tree and $G$ be a group. Then a group action of $G$ on $T$ {\it respects} $T$ for every $g\in G$,
the function $p\mapsto g.p$ is an automorphism of $T$ (which we henceforth denote by $g$). The action is {\it level-transitive} if for every
$\alpha<\height(T)$, the induced group action on $ T(\alpha)$ is transitive, meaning that for every two nodes $p,q\in T(\alpha)$, there is a $g\in G$
such that $g.p=q$.
\end{definition}

An equivalent way of saying that the action is level transitive is that the orbit $G[p]=\{g.p\st g\in G\}$ of any node $p$ in the tree under the
group action is simply the corresponding level $T(|p|)$ of the tree. It follows that for any cofinal branch through the tree, the images of this
branch under the group action yield the entire tree. In the following, we will mainly be concerned with groups of automorphisms of a tree $T$ and
their canonical actions on $T$, namely, the action $\pi.p=\pi(p)$.

\begin{definition}\rm
\label{def:respecting-groups} A set of cofinal branches $B$ through a tree $T$ {\it covers} the tree if $\bigcup B=T$, so that every node in $T$ lies
on a branch in $B$. The set $B$ {\it respects} an automorphism $\pi$ of $T$ if $B$ is closed under pointwise application of $\pi$; that is, if
whenever $b\in B$, then $\pi[b]\in B$. The set $B$ {\it respects} a group action of a group $G$ on $T$, if it respects the associated automorphisms
of every element of $G$ under the group action. If $b$ is a cofinal branch of $T$ and $g\in G$, then let $g[b]$ be the image of $b$ under the
automorphism associated to $g$. We shall write $G[b]$ for the orbit of $b$ under the group action, namely, $G[b]=\{g[b]\st g\in G\}$. Thus, $B$
respects $G$ if and only if $G[b]\sub B$, for every $b\in B$.
\end{definition}

When constructing our trees from $\diamondsuit$, we will follow the same basic strategy as in Theorem \ref{thm:DiamondImpliesRigidSuslinTree}, in
that we will use the diamond sequence to anticipate antichains or automorphisms or names of automorphisms, and so on, and then extend the tree so as
to kill off or seal these anticipated objects. By doing so, we will ensure that the $\omega_1$-tree we ultimately construct will have the desired
rigidity properties. We explain in Definition \ref{def:objects-to-seal} exactly the sense in which we will seal these various objects; one should
imagine here that we propose to extend a countable tree $T$ by adding a limit level containing nodes exactly above the branches in $B$.

\begin{definition}\rm
\label{def:objects-to-seal} Let $T$ be a tree of limit height and $B$ a set of branches that covers $T$.
\begin{enumerate}
\item $B$ \emph{seals a maximal antichain} $A$ of $T$ if for every $b\in B$, $b\cap A\neq\leer$.

\item $B$ \emph{seals a nontrivial automorphism} $\pi$ of $T$ if $B$ does not respect it. This means that there is a branch $b\in B$ such that
$\pi[b]\notin B$.

\item A function $f$ is a \emph{potential additional branch of $T$} if for some $p\in T$ it is an order preserving map $f:T_p\To T$ with:
\begin{enumerate}
  \item[(a)] $f(p)\perp p$.
  \item[(b)] $\forall\gamma<\height(T)\,\forall q\ge_T p\, \exists r\ge_T q\quad
    |f(r)|_T\ge\gamma$.
  \item[(c)] $\forall q\ge_T p\ \exists r_0,r_1\ge_Tq\quad f(r_0)\incomp
    f(r_1)$.
\end{enumerate}

\item $B$ \emph{seals a potential additional branch} if there is a $b\in B$ such that $f[b]$ determines a cofinal branch through $T$ (the closure of
$f[b]$ under $<_T$), but this branch is not in $B$.

\item A function $f$ is a \emph{potential additional automorphism} if there is a $p\in T$ such that $\dom(f)=T_p$ and:
\begin{enumerate}
  \item[(a)] For all $q\ge_Tp$, $f(q)$ is a partial automorphism of
    $T$.
  \item[(b)] $f$ is monotonic, meaning that $p\le_T q\le_T r$ implies $f(q)\sub f(r)$.
  \item[(c)] For all $q\ge_Tp$, there are $r_0$, $r_1\ge_T q$ such
    that there exists an $s\in\dom(f(r_0))\cap\dom(f(r_1))$ with the
    property that $f(r_0)(s)\neq f(r_1)(s)$.
  \item[(d)] For all $q\ge_Tp$ and all $r\in T$, there is a $q'\ge_T q$
    such that $r\in\dom(f(q'))$.
\end{enumerate}

\item $B$ \emph{seals a potential additional automorphism} $f$ if there is a $b\in B$ such that $B$ does not respect $\bigcup f[b]$. This means that
there is a branch $c\in B$ such that $c\sub\dom(\bigcup f[b])$ and $(\bigcup f[b])[c]\notin B$.

\item The notions \emph{potential additional branch of degree $n$} and \emph{potential additional automorphism of degree $n$} are defined
analogously. Thus, a potential additional branch of degree $n$ is an order preserving function $f:T_{p_0}\times\cdots\times T_{p_{n-1}}\To T$ such
that $f(\vp)\incomp p_i$ for all $i<n$ and: 
\begin{enumerate}
  \item $\forall\gamma<\height(T)\,\forall \vq\ge_{T^n} \vp\ \exists \vr\ge_{T^n} \vq\quad
    |f(\vr)|_T\ge\gamma$.
  \item $\forall \vq\ge_T\vp\ \exists \vr^0,\vr^1\ge_{T^n}\vq\quad f(\vr^0)\incomp
    f(\vr^1)$.
\end{enumerate}
\end{enumerate}
\end{definition}

\begin{lem}
\label{lem:reason-for-wanting-to-seal-objects} Let $T$ be a normal tree of limit height.
\begin{enumerate}

  \item \label{item:UBP-eq-no-pot-add-branch}
    $T$ is UBP if and only if there is no potential additional
    branch of $T$. Analogously, $T$ is
    $n$-absolutely UBP if and only if there is no potential additional
    branch of degree $n$.
  \item \label{item:no-add-automorphisms-eq-no-pot-automorphisms}
    Forcing with $T$
    adjoins no new automorphisms of $T$ if and only if there is no
    potential additional automorphism of $T$. Again, forcing with
    $T^n$ adds no new automorphism of $T$ if and only if
    there is no potential additional automorphism of
    degree $n$.
\end{enumerate}
\end{lem}

\proof Statement \ref{item:UBP-eq-no-pot-add-branch} is proved by realizing that a potential additional branch of $T$ is essentially the same as a
$T$-name for a new branch different from the generic branch. More precisely, suppose $T$ does not have the unique branch property, so that there is a
condition $p\in T$ forcing that $\tau$ is (the name of) a new cofinal branch through $T$, different from the generic branch. Since $p$ forces that
$\tau$ is different from the generic branch, we may strengthen $p$ if necessary and assume that $p\forces r\in\tau$ for some $p\perp r$. More
generally, for each $q\in T_p$, let $f(q)$ be the $T$-maximal node $r$ such that $q\forces_T$ $r\in\tau$ (such a maximal $r$ exists by the uniqueness
of limit nodes in $T$). It is now easy to see that this is a potential additional branch. By design, $f:T_p\to T$ is order preserving and $f(p)\perp
p$. For any $\gamma<\height(T)$ and any $q\in T_p$, we may extend $q$ to some $q'$ so as to decide $\tau$ beyond height $\gamma$, so $\gamma\leq
|f(q')|$. And since $\tau$ is forced to be not in $\check V$, there is no condition deciding all of it; so for every $q\in T_p$ there are extensions
$q_0$ and $q_1$ forcing specific incompatible nodes into $\tau$, so that $f(q_0)\perp f(q_1)$.

Conversely, suppose that $f:T_p\To T$ is a potential additional branch. We claim that $p$ forces that at least two new branches are added. To see
this, suppose that $b$ is a $V$-generic branch through $T$ containing $p$. Let $c$ be the closure of $f[b]$. Some simple density arguments show that
$c$ is a cofinal branch through $T$ that is not in $V$. And since $p\perp f(p)$, we know $b\neq c$. Thus, $T$ does not have the unique branch
property, establishing the first claim of \ref{item:UBP-eq-no-pot-add-branch}. The statement for arbitrary finite $n$ follows in the same way.

A similar argument establishes statement \ref{item:no-add-automorphisms-eq-no-pot-automorphisms}. Specifically, if some $p\in T$ forces that $\pi$ is
(the name of) a new automorphism of $T$, then for each $q\in T_p$ we define the corresponding partial isomorphism $f(q)$ of $T$ by:
$$f(q)(s)=t\quad\hbox{ if and only if }\quad q\forces\pi(\check s)=\check t.$$ It is easy to see that $f$ is a potential additional automorphism of
$T$. Conversely, if we have a potential additional automorphism $f$, on domain $T_p$, then we force below $p$ to add a $V$-generic branch $b$. Some
simple density arguments now establish that $\pi=\bigcup_{q\in b,q\ge p}f(q)$ is a new automorphism of $T$. The argument for general finite $n$ is
similar.\qed

Now we are ready to state our general sealing lemma.

\begin{lem}
\label{lem:basic-sealing-lemma} Let $T$ be a countable normal tree of limit height. Let $G$ be a countable group respecting $T$ and acting
level-transitively on $T$.
\begin{enumerate}
  \item
    \label{item:OrbitsAreCovers}
    If $b$ is a cofinal branch of $T$, then $G[b]$ is a countable
    set of branches covering $T$ and respecting $G$.
  \item
    \label{item:SealingAnAntichain}
    Every maximal antichain in $T$ can be sealed by a countable
    set of branches covering $T$ and respecting $G$.
    Moreover, this
    set can be chosen to be the orbit of a single cofinal branch under
    $G$.
  \item
    \label{item:SealingAnAutomorphism}
    Let $\pi$ be a nontrivial automorphism of $T$. Then $\pi$ can
    be sealed by a countable set of branches covering $T$.

    If there is a $p_0\in T$
    such that $\pi(p_0)\neq p_0$ and $G$ satisfies the
    requirement
    \claim{$(\Gamma)$}{$\forall\sigma\in G\ \forall p\ge_T p_0\ \exists
      p'\ge_T p\quad \sigma(p')\incomp\pi(p')$,}
    then there is a countable set of branches sealing $\pi$, covering
    $T$ and respecting $G$.
    Moreover, this set can be chosen to be the
    orbit of a single branch under the group action.
  \item
    \label{lem:basic-sealing-lemma-potential-branches}
    Every potential additional branch $f$ can be sealed by a countable
    set of branches covering $T$.

    If\/ $G$ satisfies condition $(\Gamma)$ at some $p_0\in T$, with $f$ replacing $\pi$, then
    there is a countable set of branches sealing $f$, covering $T$ and
    respecting $G$.
    Again, this set can be chosen to be the orbit
    of a single branch.
  \item
    \label{item:SealingAdditionalBranchesOfDegreeN}
    Every potential additional branch $f$ of
    degree $n\in\omega$ can be sealed by a countable set of branches covering $T$.
  \item
    \label{item:AvoidingCountablyManyBranches}
    For any countable set $B$ of branches through $T$, there is
    a countable set of branches $B'$, disjoint from $B$, which covers
    $T$, respects $G$,
    and is the orbit of a single branch.
\end{enumerate}
\end{lem}

\proof \ref{item:OrbitsAreCovers}.) Let $b$ be a cofinal branch of $T$, and let $B=G[b]$ be the corresponding orbit of $b$ under the action of $G$.
Clearly, $B$ is a countable set of cofinal branches respecting $G$. To see that $B$ covers $T$, suppose $s\in T$ and let $t$ be the $|s|_T^{\rm th}$
element of $b$. Since the group action is level transitive, there is $\sigma\in G$ such that $\sigma(t)=s$, and consequently, $s\in\bigcup B$.

\ref{item:SealingAnAntichain}.) Let $A$ be a maximal antichain. For each $\sigma\in G$, let:
\[ D_\sigma=\{q\in T\st\exists a\in A\quad\sigma(q)\ge_T a\}. \]
Each $D_\sigma$ is dense in $T$, viewed as a notion of forcing. To see this, for any node $p$ use the maximality of $A$ to find a $q'\ge_T \sigma(p)$
such that $q'\ge_T a$ for some $a\in A$, and then observe that $q\mdf\sigma^{-1}(q')$ is in $D_\sigma$ and $q\ge_T p$, as desired. Using this, we may
now construct a cofinal branch $b$ through $T$ meeting each of the countably many dense sets $D_\sigma$ for $\sigma\in G$.
Let $B=G[b]$. By \ref{item:OrbitsAreCovers}, $B$ is a countable cover of $T$ that respects $G$. We argue that $B$ seals $A$. For any $\sigma\in G$,
since the branch $b$ contains some $p\in D_\sigma$, we know $\sigma(p)$ is above an element of $A$, and so every $\sigma[b]\in B$ intersects $A$, as
desired.

\ref{item:SealingAnAutomorphism}.) Suppose that $\pi$ is a nontrivial automorphism of $T$. Fix any branch $b$ such that $b\neq\pi[b]$, and for each
$p\in T$ choose a branch $b_p$ such that $p\in b_p\neq\pi[b]$, which is possible because $T$ is ever-branching. The set $B=\{b_p\st p\in
T\}\cup\{b\}$, therefore, is a countable set of branches sealing $\pi$ and covering $T$, as desired. For the more general claim, suppose that a group
$G$ acts level transitively on $T$ and satisfies $(\Gamma)$ with respect to $p_0$. For $\sigma\in G$, let:
\[ D_\sigma=\{p'\in T\st\sigma(p')\neq\pi(p')\}. \]
The condition $(\Gamma)$ exactly expresses that $D_\sigma$ is dense below $p_0$. Hence, we can find a cofinal branch $b$ through $T$ such that
$p_0\in b$ and $b$ meets every $D_\sigma$ for $\sigma\in G$. As before, let $B=G[b]$. Again, $B$ is a countable cover of $T$ which respects $G$, by
\ref{item:OrbitsAreCovers}. Moreover, it seals $\pi$, because while $b\in B$, we know from $b$ meeting $D_\sigma$ that $\sigma[b]\neq\pi[b]$, and so
$\pi[b]\notin B$. So $B$ seals $\pi$, as desired.

\ref{lem:basic-sealing-lemma-potential-branches}.) The proof here is similar to that of statement \ref{item:SealingAnAutomorphism}. To seal a
potential additional branch $f$ with domain $T_p$, we will use the fact that for every $\gamma< |T|$, the following set is dense below $p$, by clause
(b) in the definition:
\[ D_\gamma=\{r\st\gamma\le |f(r)|\}. \]
We may therefore choose a cofinal branch $b$ through $T$ meeting each of these countably many sets and containing $p$. Let $c$ be the corresponding
cofinal branch determined by $f[b]$. We now continue as in the proof of \ref{item:SealingAnAutomorphism}, using $c$ in place of $\pi[b]$. In the case
where there is a group $G$ satisfying $(\Gamma)$ with respect to $f$ and $p_0$, we modify the construction in the same way as in statement
\ref{item:SealingAnAutomorphism}. That is, we construct the cofinal branch $b$ to meet not only the previous $D_\gamma$, but also to meet the dense
sets
\[ D_\sigma=\{p'\in T\st\sigma(p')\perp f(p')\}, \]
for every $\sigma\in G$. Again, by \ref{item:OrbitsAreCovers}, $B=G[b]$ is a countable set of  branches covering $T$ and respecting $G$. The set $B$
seals $f$ because $b\in B$ but the branch $c$ determined by $f[b]$ is not in $B$.

\ref{item:SealingAdditionalBranchesOfDegreeN}.) Suppose that $f$ is a potential additional branch of some degree $n\in\omega$. So the domain of $f$
is $T_{p_0}\times\cdots\times T_{p_{n-1}}$ for some $\vec p=\<p_0,\ldots,p_{n-1}>\in T^n$, where $f(\vec p)\perp p_i$ for each $i<n$. We may find
branches $b_0,\ldots,b_n\in[T]$ with $p_i\in b_i$ such that $b_0\times\cdots\times b_n$ meets each of the dense sets
\[ D_\gamma=\{\vr\in T^n\st |f(\vr)|\ge\gamma\}. \]
This is possible because each $D_\gamma$ is dense below $\vp$ in $T^n$. Let $c$ be the branch determined by $f[b_0\times\cdots\times b_{n-1}]$, and
for each $q\in T$ choose a branch $b_q\in[T]$ such that $q\in b_q\neq c$. It follows that the set $B\mdf\{b_q\st q\in T\}\cup\{b_0,\ldots,b_n\}$
covers $T$ and seals $f$, since $c\notin B$.

\ref{item:AvoidingCountablyManyBranches}.) Suppose $B$ is a countable set of branches through $T$ and $G$ is a countable group acting level
transitively on $T$. For any branch $b\in B$ and $\sigma\in G$, the set
\[ D_{\sigma,b}=\{p\in T\st \sigma(p)\notin b\} \]
is dense in $T$. Choose a cofinal branch $c$ through $T$ meeting all $D_{\sigma,b}$ for $\sigma\in G$ and $b\in B$. By \ref{item:OrbitsAreCovers},
the set $B'=G[c]$ is a countable set of branches covering $T$ and respecting $G$. It is disjoint from $B$, by construction.\qed

Let us introduce a class of automorphisms that we will use in our construction of a UBP tree that is not absolutely rigid. Our construction will
involve certain $\gamma$-trees, subtrees of ${}^{<\gamma}2$, and we will consider the following class of automorphisms, which happen all to be
restrictions of automorphisms of the full tree ${}^{\le\gamma}2$. Specifically, for any ordinal $\gamma$ and any $s\in{}^{\gamma}2$, let $\pi_s$ be
the automorphism of ${}^{\le\gamma}2$ that simply swaps the digits at the positive coordinates of $s$. More formally, define $\pi_s(t)$ by:
\[ \pi_s(t)(\alpha)=
  \left\{
    \begin{array}{l@{\qquad}l}
      1-t(\alpha) & \tx{if} \alpha\in\dom(t) \tx{and} s(\alpha)=1,\\
      t(\alpha)   & \tx{if} \alpha\in\dom(t) \tx{and} s(\alpha)=0.\\
    \end{array}
  \right.
\]
Let $\Pi^\gamma=\{\pi_s\st s\in{}^\gamma 2\}$ be the corresponding group of such automorphisms. For $a\sub\gamma$, we shall write $\chi^\gamma_a$ to
denote the characteristic function of $a$ as a subset of $\gamma$, and we shall write $\pi^\gamma_a$ for $\pi_{\chi^\gamma_a}$. Viewing $\Pi^\gamma$
as a group of automorphisms, we shall write $\kla{S}$ for the subgroup of $\Pi^\gamma$ generated by $S\sub\Pi^\gamma$. If $G$ is a subgroup of
$\Pi^\gamma$, and $p\in{}^{\le\gamma}2$, we let $G[p]=\{\sigma(p)\mid \sigma\in G\}$ be the orbit of $p$ under the canonical group action of
$\Pi^\gamma$ on ${}^{\le\gamma} 2$. We now record some simple but crucial facts about $\Pi^\gamma$.

\goodbreak
\begin{lem}
\label{lem:properties-of-nice-automorphisms}
  Let $\gamma$ be an ordinal.
  \begin{enumerate}
  \item $\pi^\gamma_a\circ\pi^\gamma_b=\pi^\gamma_{a\bigtriangleup
      b}$, for $a,b\sub\gamma$.
  \item Every element of $\Pi^\gamma$ is self-inverse.
  \item $\Pi^\gamma$ is commutative, as a group of automorphisms.
  \item \label{item:canonical-action}
    $\Pi^\gamma$ operates on itself via the group action assigning
    to each $\pi\in\Pi^\gamma$ the action
    $\pi_s\mapsto\pi_{\pi(s)}$. Denoting this action by ``.'',
    we have:
    \[ \pi_1.\pi_2=\pi_1\circ\pi_2,\]
    for $\pi_1$, $\pi_2\in\Pi^\gamma$. That is, it is the canonical group
    action of $G$ on itself.
  \item \label{item:simple-group-generation} For any set $S\sub\Pi^\gamma$,
    the generated group $\kla{S}$ consists precisely of all
    automorphisms of the form
    $s_{i_0}\circ\cdots\circ s_{i_{n-1}}$, where
    $\kla{s_\alpha\st \alpha<\kappa}$ is a fixed enumeration of
    $S$ and $i_0<\ldots<i_n<\kappa$. (The empty
    composition is taken here to be the identity.)
  \item \label{item:groups-generated-by-orbits}
    If $G$ is a subgroup of $\Pi^\gamma$ and $p\in{}^\gamma 2$, then
    \[\kla{\{\pi_q\st q\in G[p]\}}=\kla{G\cup\{\pi_p\}}.\]
\end{enumerate}
\end{lem}

\proof These are routine verifications, although we give a proof of statement \ref{item:groups-generated-by-orbits}.  Let $A\mdf\kla{\{\pi_q\st q\in
G[p]\}}$ and $B\mdf\kla{G\cup\{\pi_p\}}$. For the inclusion from left to right, assume that $\sigma\in A$. It follows by
\ref{item:simple-group-generation} that $\sigma$ is a finite composition of automorphisms of the form $\pi_{\sigma_i(p)}$, which is equal to
$\sigma_i\circ\pi_p$ by \ref{item:canonical-action}. Since the group is commutative, this composition has the form
$\sigma_0\circ\cdots\circ\sigma_{n-1}\circ(\pi_p)^n$, which is an element of $B$, as $(\pi_p)^n$ is either the identity or $\pi_p$, depending on
whether $n$ is even or odd. Conversely, every automorphism in $B$ is of the form $\sigma$ or $\sigma\circ\pi_p$, where $\sigma\in G$. In the latter
case, $\sigma\circ\pi_p=\pi_{\sigma(p)}\in A$, and in the former, $\pi_p\in A$ because the identity is in $G$ and $\sigma\circ\pi_p\in A$, so
$\sigma\in A$. Hence the composition of these is in $A$ again. But the composition is $\sigma$ since $\pi_p$ is self inverse. \qed

We are finally ready to fulfill the third requirement of Observation \ref{obs:SufficesToFind}, which we do under the assumption of $\diamondsuit$.

\begin{thm}\label{thm:abs-non-rigidity}
If\/ $\diamondsuit$ holds, then there is a Souslin tree with the unique branch property, but which is not absolutely rigid (and is absolutely
non-rigid).
\end{thm}

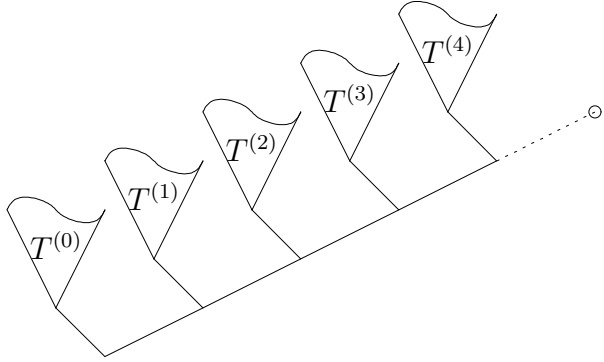
\begin{figure}[h]
\begin{center}
\setlength{\unitlength}{0.00043333in}
\begingroup\makeatletter\ifx\SetFigFont\undefined%
\gdef\SetFigFont#1#2#3#4#5{%
  \reset@font\fontsize{#1}{#2pt}%
  \fontfamily{#3}\fontseries{#4}\fontshape{#5}%
  \selectfont}%
\fi\endgroup%
{\renewcommand{\dashlinestretch}{30}
\begin{picture}(7295,4396)(0,-10)
\put(5412,3612){\makebox(0,0)[b]{\smash{{{\SetFigFont{12}{14.4}{\familydefault}{\mddefault}{\updefault}$T^{(4)}$}}}}}
\put(4212,3012){\makebox(0,0)[b]{\smash{{{\SetFigFont{12}{14.4}{\familydefault}{\mddefault}{\updefault}$T^{(3)}$}}}}}
\put(3012,2412){\makebox(0,0)[b]{\smash{{{\SetFigFont{12}{14.4}{\familydefault}{\mddefault}{\updefault}$T^{(2)}$}}}}}
\put(1812,1812){\makebox(0,0)[b]{\smash{{{\SetFigFont{12}{14.4}{\familydefault}{\mddefault}{\updefault}$T^{(1)}$}}}}}
\put(612,1212){\makebox(0,0)[b]{\smash{{{\SetFigFont{12}{14.4}{\familydefault}{\mddefault}{\updefault}$T^{(0)}$}}}}}
\path(4812,4212)(4812,4213)(4813,4215)
	(4815,4219)(4818,4225)(4822,4233)
	(4827,4242)(4833,4253)(4840,4266)
	(4849,4279)(4859,4292)(4870,4305)
	(4883,4317)(4897,4329)(4913,4340)
	(4932,4349)(4953,4357)(4977,4363)
	(5005,4367)(5037,4369)(5073,4367)
	(5112,4362)(5150,4354)(5188,4344)
	(5224,4333)(5257,4320)(5287,4306)
	(5313,4293)(5335,4279)(5355,4265)
	(5371,4252)(5386,4238)(5399,4225)
	(5412,4212)(5425,4199)(5438,4186)
	(5453,4172)(5469,4159)(5489,4145)
	(5511,4131)(5537,4118)(5567,4104)
	(5600,4091)(5636,4080)(5674,4070)
	(5712,4062)(5751,4057)(5787,4055)
	(5819,4057)(5847,4061)(5871,4067)
	(5892,4075)(5911,4084)(5927,4095)
	(5941,4107)(5954,4119)(5965,4132)
	(5975,4145)(5984,4158)(5991,4171)
	(5997,4182)(6002,4191)(6006,4199)
	(6009,4205)(6011,4209)(6012,4211)(6012,4212)
\path(12,1812)(12,1813)(13,1815)
	(15,1819)(18,1825)(22,1833)
	(27,1842)(33,1853)(40,1866)
	(49,1879)(59,1892)(70,1905)
	(83,1917)(97,1929)(113,1940)
	(132,1949)(153,1957)(177,1963)
	(205,1967)(237,1969)(273,1967)
	(312,1962)(350,1954)(388,1944)
	(424,1933)(457,1920)(487,1906)
	(513,1893)(535,1879)(555,1865)
	(571,1852)(586,1838)(599,1825)
	(612,1812)(625,1799)(638,1786)
	(653,1772)(669,1759)(689,1745)
	(711,1731)(737,1718)(767,1704)
	(800,1691)(836,1680)(874,1670)
	(912,1662)(951,1657)(987,1655)
	(1019,1657)(1047,1661)(1071,1667)
	(1092,1675)(1111,1684)(1127,1695)
	(1141,1707)(1154,1719)(1165,1732)
	(1175,1745)(1184,1758)(1191,1771)
	(1197,1782)(1202,1791)(1206,1799)
	(1209,1805)(1211,1809)(1212,1811)(1212,1812)
\path(612,612)(1212,12)(2412,612)
\dashline{45.000}(6012,2412)(7212,3012)
\path(12,1812)(612,612)(1212,1812)
\path(4812,4212)(5412,3012)(6012,4212)
\path(6012,2412)(5412,3012)
\put(7212,3012){\ellipse{150}{150}}
\path(1812,1212)(2412,612)(3612,1212)
\path(1212,2412)(1212,2413)(1213,2415)
	(1215,2419)(1218,2425)(1222,2433)
	(1227,2442)(1233,2453)(1240,2466)
	(1249,2479)(1259,2492)(1270,2505)
	(1283,2517)(1297,2529)(1313,2540)
	(1332,2549)(1353,2557)(1377,2563)
	(1405,2567)(1437,2569)(1473,2567)
	(1512,2562)(1550,2554)(1588,2544)
	(1624,2533)(1657,2520)(1687,2506)
	(1713,2493)(1735,2479)(1755,2465)
	(1771,2452)(1786,2438)(1799,2425)
	(1812,2412)(1825,2399)(1838,2386)
	(1853,2372)(1869,2359)(1889,2345)
	(1911,2331)(1937,2318)(1967,2304)
	(2000,2291)(2036,2280)(2074,2270)
	(2112,2262)(2151,2257)(2187,2255)
	(2219,2257)(2247,2261)(2271,2267)
	(2292,2275)(2311,2284)(2327,2295)
	(2341,2307)(2354,2319)(2365,2332)
	(2375,2345)(2384,2358)(2391,2371)
	(2397,2382)(2402,2391)(2406,2399)
	(2409,2405)(2411,2409)(2412,2411)(2412,2412)
\path(1212,2412)(1812,1212)(2412,2412)
\path(4212,2412)(4812,1812)(6012,2412)
\path(3612,3612)(3612,3613)(3613,3615)
	(3615,3619)(3618,3625)(3622,3633)
	(3627,3642)(3633,3653)(3640,3666)
	(3649,3679)(3659,3692)(3670,3705)
	(3683,3717)(3697,3729)(3713,3740)
	(3732,3749)(3753,3757)(3777,3763)
	(3805,3767)(3837,3769)(3873,3767)
	(3912,3762)(3950,3754)(3988,3744)
	(4024,3733)(4057,3720)(4087,3706)
	(4113,3693)(4135,3679)(4155,3665)
	(4171,3652)(4186,3638)(4199,3625)
	(4212,3612)(4225,3599)(4238,3586)
	(4253,3572)(4269,3559)(4289,3545)
	(4311,3531)(4337,3518)(4367,3504)
	(4400,3491)(4436,3480)(4474,3470)
	(4512,3462)(4551,3457)(4587,3455)
	(4619,3457)(4647,3461)(4671,3467)
	(4692,3475)(4711,3484)(4727,3495)
	(4741,3507)(4754,3519)(4765,3532)
	(4775,3545)(4784,3558)(4791,3571)
	(4797,3582)(4802,3591)(4806,3599)
	(4809,3605)(4811,3609)(4812,3611)(4812,3612)
\path(3612,3612)(4212,2412)(4812,3612)
\path(3012,1812)(3612,1212)(4812,1812)
\path(2412,3012)(2412,3013)(2413,3015)
	(2415,3019)(2418,3025)(2422,3033)
	(2427,3042)(2433,3053)(2440,3066)
	(2449,3079)(2459,3092)(2470,3105)
	(2483,3117)(2497,3129)(2513,3140)
	(2532,3149)(2553,3157)(2577,3163)
	(2605,3167)(2637,3169)(2673,3167)
	(2712,3162)(2750,3154)(2788,3144)
	(2824,3133)(2857,3120)(2887,3106)
	(2913,3093)(2935,3079)(2955,3065)
	(2971,3052)(2986,3038)(2999,3025)
	(3012,3012)(3025,2999)(3038,2986)
	(3053,2972)(3069,2959)(3089,2945)
	(3111,2931)(3137,2918)(3167,2904)
	(3200,2891)(3236,2880)(3274,2870)
	(3312,2862)(3351,2857)(3387,2855)
	(3419,2857)(3447,2861)(3471,2867)
	(3492,2875)(3511,2884)(3527,2895)
	(3541,2907)(3554,2919)(3565,2932)
	(3575,2945)(3584,2958)(3591,2971)
	(3597,2982)(3602,2991)(3606,2999)
	(3609,3005)(3611,3009)(3612,3011)(3612,3012)
\path(2412,3012)(3012,1812)(3612,3012)
\end{picture}
}
\end{center}
\caption{The trees $T^{(n)}$, glued together to form one composite tree.} \label{fig:glued-together-tree}
\end{figure}

\proof Using $\diamondsuit$, we will construct a 2-splitting Souslin tree with the unique branch property, but which is absolutely non-rigid. Our
strategy will be to construct trees $T^{(n)}$ for $n<\omega$ with the following properties:
\begin{enumerate}
 \item Each $T^{(n)}$ is a rigid 2-splitting Souslin tree.
 \item Each $T^{(n)}$ has the unique branch property.
 \item Forcing with any $T^{(n)}$ adds no cofinal branch to any other $T^{(m)}$.
 \item Forcing with $T^{(n)}$ adds a nontrivial automorphism of $T^{(n+1)}$.
\end{enumerate}
With such trees, we can build the final tree $T$ by gluing together the trees $T^{(n)}$ as in figure \ref{fig:glued-together-tree}.
Since each $T^{(n)}$ is a 2-splitting Souslin tree, the resulting glued-together tree $T$ is also a 2-splitting Souslin tree. Observe that forcing
with $T$ is equivalent to choosing some $n$ and forcing with $T^{(n)}$. Since each $T^{(n)}$ has the unique branch property and no $T^{(n)}$ adds a
branch to another $T^{(m)}$, this glued together tree $T$ therefore has the unique branch property. Since forcing with $T^{(n)}$ adds an automorphism
to $T^{(n+1)}$ and hence also to $T$, it follows that $T$ is not absolutely rigid (and even absolutely non-rigid). Thus, the glued together tree $T$
will satisfy all our desired properties, proving the theorem.

Let us begin the construction. We will construct the trees $T^{(n)}$ by simultaneous recursion on the levels. Each tree $T^{(n)}$ will consist of
binary sequences, ordered by inclusion. We will inductively maintain that $T^{(n)}|\alpha$ is a normal $\alpha$-tree, in fact a subtree of
${}^{<\alpha}2$, ordered by inclusion. In addition, for any $\gamma<\alpha<\omega_1$ and $p\in T^{(n)}(\gamma)$, we will ensure that
$\pi_p\rest(T^{(n+1)}|(\gamma+1))$ is an automorphism of $T^{(n+1)}|(\gamma+1)$, and the group generated by these automorphisms acts
level-transitively on $T^{(n+1)}|(\gamma+1)$.

Suppose that the trees $T^{(n)}|\alpha$ have been constructed below $\alpha$; we must now specify the $\alpha^{\rm th}$ levels $T^{(n)}(\alpha)$. At
successor levels, there is no choice: we give every node on the top level of the tree two immediate successors by adjoining $0$ and $1$ to the binary
sequences. This is easily seen to maintain our inductive assumptions. The interesting case occurs when $\alpha$ is a limit ordinal. In this case, our
construction will always proceed with the following template. First, we will specify a particular set of branches $B_0$ covering $T^{(0)}|\alpha$ and
particular individual branches $b_n\in[T^{(n)}|\alpha]$ for $n>0$. This information will determine the $\alpha^{\rm th}$ levels of the trees
$T^{(n)}(\alpha)$ as follows. We let $G_0$ be the group generated by $\{\pi_b\mid b\in B_0\}$. Our induction hypothesis ensures that this group acts
level transitively on $T^{(1)}|\alpha$, and so the set $B_1=G_0[b_1]$ is a covering set of branches for $T^{(1)}|\alpha$, by Lemma
\ref{lem:basic-sealing-lemma}.\ref{item:OrbitsAreCovers}. Continuing recursively, we define $G_n$ to be the group of automorphisms generated by
$B_n$, and $B_{n+1}=G_n[b_{n+1}]$. Inductively, $G_n$ acts level transitively on $T^{(n+1)}|\alpha$, and so $B_{n+1}$ covers $T^{(n+1)}|\alpha$. We
now extend the trees to level $\alpha$ by defining $T^{(n)}(\alpha)=B_n$, conflating branches through ${}^{<\alpha}2$ with elements of ${}^\alpha 2$.
Since $B_{n+1}$ is the orbit of a single branch under $G_n$, it follows that $G_n$, which is the group generated by the elements of
$T^{(n)}(\alpha)=B_n$, acts level transitively on $T^{(n+1)}(\alpha)$. Therefore, as long as our construction follows this pattern, we will preserve
our induction hypotheses.

We therefore remain relatively free to choose the initial set of branches $B_0$ covering $T^{(0)}|\alpha$ and the individual branches $b_n\in
[T^{(n)}|\alpha]$ for $n>0$. We will do so in a way that will ensure that the trees $T^{(n)}$ are Souslin, that they each have the unique branch
property and more generally, that forcing with $T^{(n)}$ will not add branches through any other $T^{(m)}$. In order to accomplish this, we will
anticipate via $\diamondsuit$ the relevant potential additional branches and then seal them. So let us begin the detailed construction. Suppose that
$\<D_\alpha\mid\alpha<\omega_1>$ is a $\diamondsuit$ sequence. We assume that the trees $T^{(n)}|\alpha$ are defined up to the limit ordinal
$\alpha$; we must now specify the branches $B_0\of[T^{(0)}|\alpha]$ and the individual branches $b_n\in [T^{(n)}|\alpha]$.

Case 0. We act first to ensure that the trees $T^{(n)}$ are Souslin. Suppose that $D_\alpha$ codes a set of the form $\<0,n,A>$, where $A$ is a
maximal antichain in $T^{(n)}|\alpha$. If $n=0$, then we may choose a countable cover $B_0$ of $T^{(0)}|\alpha$ that seals $A$, and choose $b_n\in
[T^{(n)}|\alpha]$ for $n>0$ arbitrarily. Since our construction template leads to $T^{(0)}(\alpha)=B_0$, this will seal the antichain $A$ below level
$\alpha$. So assume $n>0$. In this case, we choose an arbitrary cover $B_0$ of $T^{(0)}|\alpha$, and arbitrary branches $b_1,\ldots, b_{n-1}$ through
$(T^{(1)}|\alpha), \ldots, (T^{(n-1)}|\alpha)$, respectively. Let $B_m$ and $G_m$ be the resulting covering sets of branches and the corresponding
groups in the construction template, for $m<n$. In particular, the group $G_{n-1}$ acts level transitively on $T^{(n)}|\alpha$. By Lemma
\ref{lem:basic-sealing-lemma}.\ref{item:SealingAnAntichain}, there is a branch $b_n\in[T^{(n)}|\alpha]$ such that the corresponding set of branches
$B_n=G_n[b_n]$ covers $T^{(n)}|\alpha$, respects $G_n$ and seals $A$. The remaining branches $b_k\in[T^{(k)}|\alpha]$ for $k>n$ may be chosen
arbitrarily, and we correspondingly define the $\alpha^{\rm th}$ level of the trees according to the construction template. This procedure will
ensure that the ultimate trees $T^{(n)}$ we construct will all be Souslin, because if $A\of[T^{(n)}]$ is a maximal antichain, then there will be a
stationary set of $\alpha$ for which $D_\alpha$ codes the triple $\<0,n,A\intersect T^{(n)}|\alpha>$, at which point we will seal $A$ below $\alpha$.
Therefore, it must be that $A=A\intersect T^{(n)}|\alpha$, and so the antichain was countable.

Case 1. Next, we act to ensure that the initial tree $T^{(0)}$ has the unique branch property. Suppose that $D_\alpha$ codes $\<1,f>$, where $f$ is a
potential additional branch for $T^{(0)}|\alpha$. By Lemma \ref{lem:basic-sealing-lemma}, we may choose a countable set of branches $B_0$ that covers
$T^{(0)}|\alpha$ and seals $f$. Choose branches $b_n\in[T^{(n)}|\alpha]$ for $n>0$ arbitrarily and follow the construction template to define the
$\alpha^{\rm th}$ level of the trees. Since this results in $T^{(0)}(\alpha)=B_0$, we have sealed $f$. It follows that the ultimate tree $T^{(0)}$ we
construct will have the unique branch property, for if it did not, then there would be a potential additional branch $f$ for $T^{(0)}$, and for a
stationary set of $\alpha$ we would have that $f\restrict (T^{(0)}|\alpha)$ is a potential additional branch for $T^{(0)}|\alpha$ and
$D_\alpha=\<1,f\restrict (T^{(0)}|\alpha)>$. At such a stage, we will have sealed $f\restrict T^{(0)}|\alpha$, meaning that there is a branch $b\in
T^{(0)}(\alpha)$ such that $f[b]$ is cofinal in $T^{(0)}|\alpha$ but has no upper bound in $T^{(0)}(\alpha)$. This contradicts that $f$ was a
potential additional branch on the entire tree $T^{(0)}$. So there can be no such potential additional branch, and the resulting tree $T^{(0)}$ will
have the unique branch property.

Case 2. Next, we act to ensure that forcing with an earlier tree $T^{(m)}$ will not add a branch to a later tree $T^{(n)}$, for $m<n$. Suppose that
$D_\alpha$ codes $\<2,m,n,f>$, where $m<n$ and $f$ is a potential additional branch for $(T^{(m)}|\alpha)\sqcup (T^{(n)}|\alpha)$ with $\dom(f)\of
T^{(m)}|\alpha$ and $\ran(f)\of T^{(n)}|\alpha$. Choose a countable cover $B_0$ of $T^{(0)}|\alpha$, taking care that if $m=0$, then there is a
$b_0\in B_0$ such that $f[b_0]$ is cofinal in $T^{(n)}|\alpha$. One can find such a branch by meeting countably many dense sets, since it is dense in
$T^{(m)}|\alpha$ that the values of $f(q)$ grow unbounded in $\alpha$. Next, choose cofinal branches $b_1,\ldots,b_{n-1}$ through
$T^{(1)}|\alpha,\ldots, T^{(n-1)}|\alpha$, taking care that the branch $b_m$ is chosen so that $f[b_m]$ is cofinal in $T^{(n)}|\alpha$. Let
$B_1\ldots,B_{n-1}$ be the resulting covering sets of branches, with the associated groups $G_1,\ldots,G_{n-1}$. We choose the next branch $b_n$ with
a bit more care. Namely, let $c$ be the cofinal branch of $T^{(n)}|\alpha$ generated by $f[b_m]$, and by
\ref{lem:basic-sealing-lemma}.\ref{item:AvoidingCountablyManyBranches}, pick $b_n$ in such a way that $G_{n-1}[b_n]\cap\{c\}=\leer$, thereby avoiding
the branch $c$. It follows that the induced set of branches $B_n=G_{n-1}[b_n]$ seals $f$, since $b_m$ is in $B_m$ and is hence extended, while the
corresponding branch $c$, generated by $f[b_m]$, is not. Finally, we complete the construction by choosing the remaining branches $b_k$ for $k>n$
arbitrarily, and building the corresponding $\alpha^{\rm th}$ levels of the trees according to the construction template. This procedure will ensure
that forcing with the ultimate tree $T^{(m)}$ will not add a branch to $T^{(n)}$, for if there were a $T^{(m)}$-name $\tau$ for such a branch through
$T^{(n)}$, then there would be a potential additional branch $f$ with $\dom(f)\of T^{(m)}$ and $\ran(f)\of T^{(n)}$, and for a stationary set of
$\alpha$ the restriction $f\restrict (T^{(m)}|\alpha)$ would be a potential additional branch of the kind we considered in this case and the set
$D_\alpha$ would code $\<2,m,n,f\restrict (T^{(m)}|\alpha)>$. At such a stage $\alpha$, we would have sealed $f$ below $\alpha$ by adding the branch
$b_m$ to $T^{(m)}(\alpha)$ but not extending $f[b_m]$ to any node in $T^{(n)}(\alpha)$, contradicting the fact that $f$ was a potential additional
branch.

Case 3. Next, we act to ensure that forcing with a later tree $T^{(m)}$ will not add a branch to an earlier tree $T^{(n)}$, for $n<m$. Suppose that
$D_\alpha$ codes $\<3,m,n,f>$, where $n<m$ and $f$ is a potential additional branch for $(T^{(m)}|\alpha)\sqcup (T^{(n)}|\alpha)$ with
$\dom(f)=(T^{(m)}|\alpha)_\tp$ and $\ran(f)\of T^{(n)}|\alpha$. Choose any set of branches $B_0$ covering $T^{(0)}|\alpha$, and any branches
$b_1,\ldots,b_{m-1}$ cofinal in $(T^{(1)}|\alpha),\ldots,(T^{(m-1)}|\alpha)$, respectively. Let $B_1,\dots,B_{m-1}$ and $G_0,\ldots,G_{m-1}$ be the
corresponding sets of branches and automorphism groups resulting from our construction template. We choose the next branch $b_m\in[T^{(m)}|\alpha]$
in such a way so as to seal $f$. For any $\sigma\in G_{m-1}$, $c\in B_n$ and $\gamma<\alpha$, let
$$D_{\sigma,c,\gamma} = \Bigl\{\ q\in (T^{(m)}|\alpha)\quad\Big|\qquad
 \lower8pt\vbox{\hbox{$\gamma\le|q|$, and if $\sigma(q)\in\dom(f)$,}\smallskip
       \hbox{then $f(\sigma(q))\perp c$ and $|f(\sigma(q))|\ge\gamma$}}\quad\Bigr\}.$$
We argue that this set is dense in $T^{(m)}|\balpha$. Given any $\bp\in T^{(m)}|\alpha$, choose an extension $p$ of $\bp$ such that $|p|\ge\gamma$.
If $p\in D_{\sigma,c,\gamma}$, we're done. Otherwise, since $f$ is a potential additional branch and $\sigma(p)\in\dom(f)$, there is a node
$p'\ge\sigma(p)$, such that $|f(p')|\ge\gamma$. Further, there are extensions $r_0$ and $r_1$ of $\sigma(p')$ in $T^{(m)}$ such that $f(r_0)\incomp
f(r_1)$. We may assume $f(r_0)\perp c$. Let $q=\sigma^{-1}(r_0)$ (which is the same as $\sigma(r_0)$), and observe that $q\in D_{\sigma,c,\gamma}$
and $q>p\ge\bp$. So $D_{\sigma,c,\gamma}$ is indeed dense.

Continuing with the construction, we now choose $b_m$ through $T^{(m)}|\alpha$ so as to meet every $D_{\sigma,c,\gamma}$. Choose $b_k$ cofinal in
$T^{(k)}|\alpha$ for $k>m$ arbitrarily, and carry out the construction template to define the $\alpha^{\rm th}$ level of all of the trees. We will
argue that $f$ is sealed by the resulting $T^{(m)}(\alpha)=B_m=\{\sigma[b]\st\sigma\in G_{m-1}\}$. Suppose that $b\in B_m$ and $\tp\le b$. The branch
$b$ has the form $b=\sigma[b_m]$ for some $\sigma\in G_{m-1}$. Since $b_m$ meets every $D_{\sigma,c,\gamma}$ for $c\in B_n$, it follows that
$f[\sigma(b_m)]\perp c$. Thus, $f[b]$ is not extended by any branch in $B_n$, and so we have sealed $f$. This procedure in our construction therefore
ensures that forcing with the ultimate tree $T^{(m)}$ will not add a cofinal branch to $T^{(n)}$, because if there were a $T^{(m)}$-name for such a
branch, then there would be a potential additional branch function $f$, which would be anticipated and sealed at some stage $\alpha$.

Case 4. Lastly, we act to ensure that every tree $T^{(n)}$, for $n>0$, has the unique branch property. This case is more complicated than the earlier
cases. Suppose that $D_\alpha$ codes $\<4,n,f>$, where $n>0$ and $f$ is a potential additional branch for $T^{(n)}|\alpha$. We will choose a set of
branches $B_0=\{d_k\mid k<\omega\}$ covering $T^{(0)}|\alpha$ and individual branches $b_m\in[T^{(m)}|\alpha]$ for $m>0$ in such a way that, after
following the construction template, the resulting tree $T^{(n)}(\alpha)$ seals $f$. Choosing these branches amounts to choosing a filter in the
following poset $\P$, defined with finite support in each factor:
 $$\P=(T^{(0)}|\alpha)^\omega\times\prod_{m>0} T^{(m)}|\alpha.$$
Such a filter $H\of\P$ determines $B_0=\{d_k\mid k<\omega\}$ and $b_m$ for $m>0$ by simply projecting onto the appropriate coordinates, so that $d_k$
is the projection of $H$ onto the $k^{\rm th}$ coordinate of the first factor of $\P$ and $b_m$ is the projection onto the $m^{\rm th}$ coordinate of
the second factor. We will construct $H$ to have the desired properties by meeting countably many dense sets in $\P$. First, it is easy to see that
with countably many dense sets in $\P$ we can ensure that the resulting $B_0$ is a covering set of cofinal branches through $T^{(0)}|\alpha$ and that
each $b_m$ is a cofinal branch in $T^{(m)}|\alpha$. The reader can verify that an additional list of dense sets will ensure that $f[b_n]$ is cofinal
in $T^{(n)}|\alpha$. In order to seal $f$, what we want to do is ensure that we do not add a branch to $B_n$ extending $f[b_n]$. Of course, $B_n$ is
determined by $b_n$ and $G_{n-1}$, which is determined by the earlier choices of $b_m$ and so on, and ultimately also by $B_0$.

Let us briefly analyze in greater detail the groups $G_m$ arising in the construction template. We claim that each $G_m$ is the group generated by
$G_0$ and the individual automorphisms $\pi_{b_1},\ldots,\pi_{b_m}$. We prove this by induction on $m$. For the anchor case $m=0$, there is nothing
to show. Assume inductively that the claim holds for $m$, and consider $G_{m+1}=\{\pi_b\mid b\in B_{m+1}\}=\{\pi_b\mid b\in G_m[b_{m+1}]\}$. By Lemma
\ref{lem:properties-of-nice-automorphisms}, item \ref{item:groups-generated-by-orbits}, it follows that $G_{m+1}=\kla{G_m\cup\{\pi_{b_{m+1}}\}}$. By
induction, $G_m$ is generated by $G_0$ and $\pi_{b_1},\ldots,\pi_{b_m}$, so it now follows that $G_{m+1}$ is generated by $G_0$ and
$\pi_{b_1},\ldots,\pi_{b_m},\pi_{b_{m+1}}$, as we claimed. It follows that $T^{(n)}(\alpha)=B_n$ will be the same as
$\kla{G_0\cup\{\pi_{b_1},\ldots,\pi_{b_{n-1}}\}}[b_n]$.

We now continue with the construction. We want to ensure that $f[b_n]$ is not extended to any branch in $B_n$. Thus, it will suffice to ensure that
$f[b_n]$ is incompatible with $\sigma[b_n]$ for every $\sigma\in\<G_0\cup\{\pi_{b_1},\ldots,\pi_{b_{n-1}}\}>$. Any such $\sigma$ has the form
$\sigma=(\pi_{d_{i_1}}\circ\cdots\circ\pi_{d_{i_u}})\circ(\pi_{b_{j_1}}\circ\cdots\circ\pi_{b_{j_v}})$, where $a=\{i_1,\ldots,i_u\}\of\omega$ and
$b=\{j_1,\ldots,j_v\}\of\{1,\ldots,n-1\}$ (allowing $a$ or $b$ to be empty). The conditions in $\P$ provide partial information about the
corresponding branches $\vec d$ and $\vec b$, and we will choose the filter $H$ in such a way so as to ensure that $\sigma[b_n]\perp f[b_n]$. For
notational convenience let us regard elements $p\in\P$ as having the form $p=(\vec p, \vec q)$, where $\vec p=\<p_0,p_1,\ldots>$ and $\vec
q=\<q_1,q_2,\ldots>$, with $p_i\in(T^{(0)}|\alpha)$ providing partial information about $d_i$ and $q_j\in (T^{(j)}|\alpha)$ providing partial
information about $b_j$. Note that all but finitely many $p_i$'s and $q_j$'s will be the root node of the corresponding tree, since the product $\P$
used finite support.

For any such condition $p=(\vec p,\vec q)$ and any $a=\{i_1,\ldots,i_u\}\of\omega$ and $b=\{j_1,\ldots,j_v\}\of\{1,\ldots,n-1\}$, let
$\pi^{a,b}_{\vec p,\vec q}$ be the corresponding partial partial automorphism
\[\pi^{a,b}_{\vec p,\vec q}\ \ \mdf\ \ (\pi_{p_{i_1}}\circ\cdots\circ\pi_{p_{i_u}})\circ(\pi_{q_{j_1}}\circ\cdots\circ\pi_{q_{j_v}}).\]
Such a partial automorphism is exactly what will grow into a $\sigma\in G_{n-1}$ as we explained above. In the following, we shall allow ourselves to
apply $\pi_s$ to $t$ even if $|t|>|s|$, by defining $\pi_s(t)=\pi_s(t\rest|s|)$. It follows that if $s\sub s'$, then $\pi_s(t)\sub\pi_{s'}(t)$. Let
$\dom(f)=(T^{(n)}|\alpha)_\tp$. We claim now that for any $a=\{i_1,\ldots,i_u\}\of\omega$ and $b=\{j_1,\ldots,j_v\}\of\{1,\ldots,n-1\}$, the
following set $D_{a,b}$ is dense in $\P$.
\[ D_{a,b}\mdf\bigl\{\ (\vec p,\vec q)\in\P\ \st\
q_n\perp\tp\text{ \ or \ }\pi^{a,b}_{\vec p,\vec q}(q_n)\incomp f(q_n)\ \bigr\}\] If $a$ and $b$ are both empty, then the claim follows from the fact
that $\tp\perp f(\tp)$, as $f$ is a potential additional branch. If $a$ or $b$ is nonempty, then any change to any $p_i$ or $q_j$ appearing in
$a\union b$ at an ordinal below $|q_n|$ will cause a change in $\pi^{a,b}_{\vec p,\vec q}(q_n)$ at the corresponding ordinal. Thus, some of these
extensions will exhibit incompatibility with $f(q_n)$, putting them in $D_{a,b}$. More specifically, suppose $p=(\vec p,\vec q)\in\P$ is given. If
$q_n\perp\tp$, then we're done; so assume the contrary. Suppose $a$ is nonempty, and fix some $i\in a$ (the case where $b$ is nonempty is similar).
Let $\delta=|p_i|$. By extending $p$ to a stronger condition, we may assume that $p_i$ is the shortest sequence, that is, the lowest node, appearing
in $p$, so that all other $p_k$ and $q_l$ are defined at $\delta$. Similarly, by extending $q_n$ inside $\vec q$, we may assume that $|f(q_n)|$ is
also larger than $\delta$. Now build two extensions $r$ and $s$ of $p$ by extending $p_i$ to add either $0$ or $1$. That is, $r=(\vec r,\vec q)$ and
$s=(\vec s,\vec q)$, where $r_i=p_i\verl\<0>$ and $s_i=p_i\verl\<1>$, and otherwise $r_k=s_k=p_k$ for $k\neq i$. Because we made just one change, at
$\delta$, it follows that $\pi^{a,b}_r(q_n)$ and $\pi^{a,b}_s(q_n)$ differ at $\delta$. Thus, they cannot both be compatible with $f(q_n)$, and so
$r$ or $s$ must be in $D_{a,b}$. So $D_{a,b}$ is dense, as we claimed.

We now complete the argument in this case. Let $H$ be a filter in $\P$ meeting all the dense sets we have mentioned, such that $\tp$ is in the
projection of $H$ onto the $n^{\rm th}$ coordinate of the second factor of $\P$ (so that the fact that it meets $D_{a,b}$ is meaningful). Let
$B_0=\{d_k\mid k<\omega\}$ and $b_m$ for $m>0$ be the resulting branches, obtained by projecting $H$. By construction, $B_0$ covers $T^{(0)}|\alpha$
and $b_m$ is cofinal in $T^{(m)}|\alpha$. Furthermore, any $\sigma\in G_{n-1}$ is determined, as we explained, by some $a$ and $b$ as above (so that
$\sigma=\pi^{a,b}_{\vec d,\vec b}$), and since $H$ meets $D_{a,b}$ we ensured that $\sigma[b_n]\perp f[b_n]$. Thus, we have $b_n\in B_n$ but $f[b_n]$
is not extended by any element of $B_n$. Therefore, our construction template ensures that $f$ is sealed below $\alpha$. It follows that the ultimate
tree $T^{(n)}$ we construct will have the unique branch property, since as in the other cases, any putative potential additional branch reflects to a
stage $\alpha$, where it was sealed.

We have therefore constructed the trees $T^{(n)}$ to have the four features we claimed at the beginning of this proof. So the proof of Theorem
\ref{thm:abs-non-rigidity} is now complete. \qed

With Theorems \ref{thm:AbsolutelyRigidNotTotallyRigid}, \ref{thm:AbsolutelyTotallyRigidNotUBP} and \ref{thm:abs-non-rigidity}, we have now fulfilled
the three requirements of Observation \ref{obs:SufficesToFind}, using $\diamondsuit$. Consequently, our main result, Main Theorem
\ref{thm:MainTheorem}, is now proved.

\section{A larger context of rigidity}\label{Section:LargerContext}

The implication diagram of Figure \ref{fig:ImplicationDiagram} is part of a larger implication diagram, pictured below, featuring the other rigidity
notions we have considered. The diagram continues to the right by considering $n$-absolute forms of the rigidity notions, leading up to
${<}\omega$-absolute rigidity notions or more.
{\small \begin{figure}[hc]
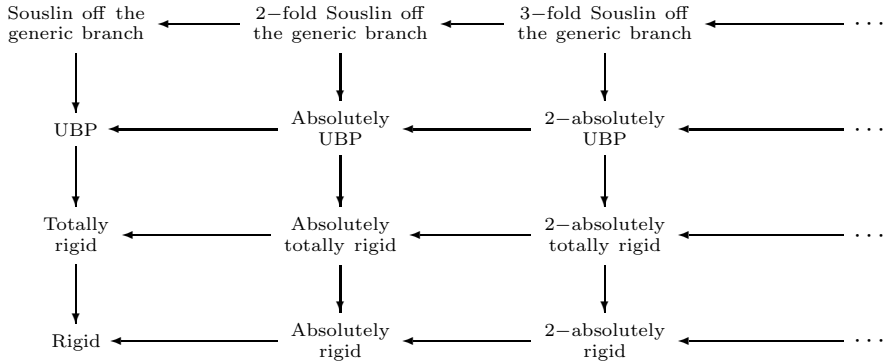

\begin{diagram}[width=5em,height=2em,objectstyle=\rm]
 {Souslin\ off\ the\atop generic\ branch} &   \lTo        &  {2-fold\ Souslin\ off\atop the\ generic\ branch}& \lTo     & {3-fold\ Souslin\ off\atop the\ generic\  branch}   & \lTo & \cdots \\
 \dTo                               &               &  \dTo                                         &          & \dTo                                     &      &    \\
 {\scriptstyle UBP}                 &   \lTo        &  {Absolutely\atop UBP}                        & \lTo     & {2-absolutely\atop UBP}                  & \lTo & \cdots  \\
 \dTo                               &               &  \dTo                                         &          & \dTo                                              \\
 {Totally\atop rigid}               &   \lTo        &  {Absolutely\atop totally\ rigid}             & \lTo     & {2-absolutely\atop totally\ rigid}        & \lTo & \cdots \\
  \dTo                              &               &  \dTo                                         &          & \dTo                                         \\
 {\scriptstyle Rigid}               &  \lTo         &  {Absolutely\atop rigid}                      & \lTo     & {2-absolutely\atop rigid}                & \lTo & \cdots \\
\end{diagram}
\caption{Larger implication diagram}\label{fig:LargerImplicationDiagram}
\end{figure}
}

\begin{question} Is the implication diagram of Figure \ref{fig:LargerImplicationDiagram} complete?
\end{question}

It seems quite possible to modify our $\diamondsuit$ constructions so as to ensure that the resulting tree is Souslin off the generic branch, by
simply anticipating and sealing the potential antichains. Conversely, to construct a tree that is UBP but not Souslin off the generic branch, one
idea would be to make a Souslin tree such that adding a branch through it creates a specializing function for some other part of the tree. Beyond
this, it seems difficult to separate Souslin off the generic branch from absolutely UBP. We leave such problems for the future.

\section{An application to the automorphism tower problem}\label{Section:AutoTowers}

In this final section, we present an application of our methods to a question involving the automorphism tower problem in group theory, an
application which was the original motivation of our investigation. The issue is that while the main theorem of
\cite{HamkinsThomas2000:ChangingHeights} showed that it is consistent with \ZFC\ that there is a group whose automorphism tower is highly malleable
by forcing, it was not known whether such a group exists in all models of set theory or, for example, in the constructible universe $L$. The specific
group $G$ constructed in \cite{HamkinsThomas2000:ChangingHeights} was obtained by first forcing to add several $V$-generic Souslin trees and then
constructing $G$ from a certain subgroup of the automorphism group of a certain graph built from these trees. The strong rigidity properties of the
generic Souslin trees led to the desired malleability of the automorphism tower of this group. Here, we use $\diamondsuit$ to construct such Souslin
trees with the desired rigidity properties. Consequently, in $L$ there are such groups whose automorphism towers are highly malleable by forcing.

We quickly review the automorphism tower construction in group theory. If $G$ is a group, then so is the automorphism group $\Aut(G)$, the set of
isomorphisms of $G$ to itself, and there is a natural homomorphism of $G$ into $\Aut(G)$ arising from conjugation. Specifically, every element $g\in
G$ maps to the corresponding inner automorphism $i_g:h\mapsto ghg^{-1}$. The automorphism tower of $G$ is obtained by iterating this process
transfinitely.
$$G_0\to G_1\to\cdots\to G_\alpha\to G_{\alpha+1}\to\cdots$$
One begins with $G_0=G$ and uses the canonical homomorphism of $G_\alpha$ into $G_{\alpha+1}=\Aut(G_\alpha)$ at successor steps. At any limit stage,
$G_\lambda$ is the direct limit of the previous groups $G_\alpha$, for $\alpha<\lambda$, with respect to these natural homomorphisms. The
automorphism tower {\it terminates} if it eventually reaches a fixed point, that is, if for some $\alpha$, the canonical map from $G_\alpha$ to
$G_{\alpha+1}$ is an isomorphism. This occurs if $G_\alpha$ is a {\it complete} group, a centerless group having only inner automorphism. The {\it
height} of the tower is the least $\alpha$ for which this occurs. If $G$ is centerless, then one can show that all the groups $G_\alpha$ in the tower
are centerless, and consequently all the maps $G_\alpha\to G_{\alpha+1}$ are injective. In this case, therefore, one can identify each group with its
image and view the tower as building up larger and larger groups, with direct limits corresponding simply to unions. Wielandt \cite{Wielandt} proved
the classical result that the automorphism tower of any centerless finite group terminates in finitely many steps. Later, various authors
\cite{RaeRoseblade}, \cite{Hulse} proved that larger classes of centerless groups had terminating automorphism towers. Simon Thomas
\cite{Thomas:AutomorphismTowerProblem}, \cite{Thomas:AutomorphismTowerProblemII} solved the automorphism tower problem for centerless groups by
proving that every centerless group has a terminating automorphism tower. Building on this, the second author \cite{Hamkins98:EveryGroup} proved that
every group has a terminating automorphism tower. The main theorem of \cite{HamkinsThomas2000:ChangingHeights} revealed that the automorphism tower
problem has what might be termed a set theoretic essence, namely, the fact that the automorphism tower of a group can be sensitive to the model of
set theory in which it is computed. We state a special case of this phenomenon here:
\begin{thm}[\cite{HamkinsThomas2000:ChangingHeights}]\label{Theorem.HT98omega}
For any $n<\omega$, there is a forcing extension with a group $G$, whose automorphism tower has height exactly $n$, but for any nonzero $m<\omega$,
there is a further (c.c.c.) forcing extension in which the automorphism tower of the very same group $G$ has height $m$.
\end{thm}
Thus, the automorphism tower of $G$ is sensitive to the set theoretic background, and even the height of the automorphism tower of $G$ can be
precisely controlled, becoming shorter or taller in various forcing extensions, as desired. The actual theorem proved in
\cite{HamkinsThomas2000:ChangingHeights} is stronger than stated above and is obtained by replacing $\omega$ with an arbitrary ordinal $\lambda$;
indeed, \cite{HamkinsThomas2000:ChangingHeights} shows that it is consistent to have a proper class of such groups for all ordinals $\lambda$ at
once. The open question with these theorems was whether one must force to add the group $G$. Perhaps one can prove in \ZFC\ that there are such
groups? Are there such groups in the constructible universe $L$? We answer at least this latter question by proving that the combinatorial principle
$\diamondsuit$, which holds in $L$, is sufficient to construct the groups.

\begin{thm}\label{Theorem.DiamondImpliesGroups}
Assume $\diamondsuit$ holds. Then for every $n<\omega$ there is a group $G$, whose automorphism tower has height $n$, but for any nonzero $m<\omega$
there is a (c.c.c.) forcing extension in which the automorphism tower of $G$ has height $m$.
\end{thm}

Let us state without much explanation that the main group-theoretic construction of \cite{HamkinsThomas2000:ChangingHeights} shows that Theorem
\ref{Theorem.DiamondImpliesGroups} is a consequence of the following combinatorial Theorem \ref{Theorem.DiamondImpliesTrees}, which we will prove.
The group $G$ of Theorem \ref{Theorem.DiamondImpliesGroups} is closely connected with a certain subgroup of the automorphism group of a graph
containing multiple copies of the trees $T^n$ of Theorem \ref{Theorem.DiamondImpliesTrees}. The crucial property 3 of Theorem
\ref{Theorem.DiamondImpliesTrees} allows the automorphism group of this graph to be precisely controlled by forcing. For the details of this
construction, we refer the readers to \cite{HamkinsThomas2000:ChangingHeights}, and to the survey article \cite{Hamkins2001:HowTall?}, which contains
a gentle overview of it. Thomas' forthcoming monograph \cite{Thomas:AutomorphismTowerProblemBook} includes extensive discussion of this construction
and many other aspects of the automorphism tower problem.

\begin{thm}\label{Theorem.DiamondImpliesTrees}
Assume $\diamondsuit$ holds. Then there is an infinite sequence of Souslin trees $T^n$ with the following properties:
\begin{enumerate}
 \item Each $T^n$ is a rigid Souslin tree.
 \item The trees $T^n$ are pairwise non-isomorphic.
 \item For any $m<\omega$, there is a c.c.c., countably distributive forcing extension preserving the rigidity of all of the trees, in which $T^0\cong\cdots\cong T^m$, but
         $T^i\not\cong T^j$ if $m\leq i<j$.
\end{enumerate}
\end{thm}

\proof Suppose that $\vec D=\<D_\alpha\mid\alpha<\omega_1>$ is a $\diamondsuit$-sequence. We will define the trees $T^n$ by recursively defining
their levels $T^n(\alpha)$. Simultaneously, we will recursively define certain controller trees $T^{n,m}$, for every $n<m<\omega$, which will be used
when it is desired to force the existence of an isomorphism from $T^n$ to $T^m$. All of these trees will be rigid Souslin trees. The difficulty will
be to ensure that forcing with a finite product of the controller trees, in order to force $T^0\cong\cdots\cong T^m$ as in statement 3 of Theorem
\ref{Theorem.DiamondImpliesTrees}, will not create unwanted automorphisms of any $T^n$ or unwanted isomorphisms from $T^i$ to $T^j$ for $m\leq i<j$.

We will construct the trees $T^n$ and $T^{n,m}$ so as to inductively maintain the following properties, where $\pi_s$ is the automorphism swapping
the values at the positive coordinates of $s$, as in Lemma \ref{lem:properties-of-nice-automorphisms}.
\begin{enumerate}
 \item Each $T^n|\alpha$ and $T^{n,m}|\alpha$ is a normal $\alpha$-tree, a subtree of ${}^{<\alpha}2$.
 \item If $s\in T^{n,m}(\alpha)$, then $\pi_s\restrict (T^n|\alpha)$ is an isomorphism of $T^n|\alpha$ with $T^m|\alpha$.
 \item These actions are level transitive in the following sense: if $p\in T^n(\alpha)$ and $q\in T^m(\alpha)$, then there is a finite list
     $n_0,n_1,\ldots,n_k$
     leading from $n_0=n$ to $n_k=m$, and $s_i\in
     T^{n_i,n_{i+1}}(\alpha)$ such that
     $\pi_{s_{k-1}}\compose\cdots\compose\pi_{s_0}(p)=q$.
     Here, we let $T^{m,n}\mdf T^{n,m}$ in case $m>n$. Note that this
     makes sense, since every $\pi_s$, coming from $s\in T^{n,m}$, is
     also an isomorphism from $T^m$ to $T^n$, since it is
     self-inverse.
\end{enumerate}
All the trees begin, of course, with the empty root node $\<>$, and at successor stages we extend every node on the $\alpha^{\rm th}$ level with its
two immediate successors in ${}^{\alpha}2$. It is easy to see that this maintains our inductive assumptions. If $\lambda$ is a limit ordinal and the
trees are defined at levels below $\lambda$, then we will always define the $\lambda^{\rm th}$ level by specifying for each pair $n,m$ a covering set
of branches $B^{n,m}\of [T^{n,m}|\lambda]$ and, for one specific $n_0$, a path $b\in [T^{n_0}|\lambda]$. We will then define the $\lambda^{\rm th}$
level of the controller trees to be $T^{n,m}(\lambda)=B^{n,m}$ and the $\lambda^{\rm th}$ level $T^n(\lambda)$ to consist of all images $\pi_{\vec
s}(b)$ under the resulting compositions $\pi_{\vec s}=\pi_{s_{n-1}}\compose\cdots\compose\pi_{s_0}$ leading from $T^{n_0}|\lambda$ to $T^n|\lambda$,
with $s_i\in B^{n_i,n_{i+1}}$, as in our inductive assumption 3 above.

We claim that as long as our recursive construction follows this pattern, then we will maintain our inductive assumptions. To see this, observe first
that the tree $T^{n,m}|(\lambda+1)$ will be normal, because $B^{n,m}$ covers $T^{n,m}|\lambda$. The tree $T^n|(\lambda+1)$ will be normal because the
level-transitive action of the automorphisms ensures that the images of $b$ cover every tree. Every $s\in T^{n,m}(\lambda)$ will provide an
isomorphism $\pi_s$ of $T^n|(\lambda+1)$ to $T^m|(\lambda+1)$ since we explicitly closed $T^m(\lambda)$ under the images of such $\pi_s$. And
finally, the action will still be level-transitive, because we define the level $T^n(\lambda)$ to be the images of the single node $b$ under all the
compositions of isomorphisms. So, given $\pi_{\vs}(b)$ and $\pi_{\vt}(b)$, we can let ${\vec s}^{-1}$ be the path $\vec s$ reversed, and let $\vec
u=\vs^{-1}\verl\vt$ be the composition of the paths, leading from where $\vs$ led to where $\vt$ leads. Then, $\pi_{\vt}(b)=\pi_{\vec
u}(\pi_{\vs}(b))$, as wished. So as in Theorem \ref{thm:abs-non-rigidity} we are relatively free to choose the branches $B^{n,m}\of[T^{n,m}|\lambda]$
and the branch $b\in[T^{n_0}|\lambda]$ so as to ensure that all the trees are Souslin, rigid, mutually rigid, and so on. For this, we will use the
diamond sequence $\vec D$ as in our earlier constructions to anticipate and then seal any unwanted antichains or automorphisms or potential
automorphisms that might arise. So let us now explain the particular details of how we choose $B^{n,m}$, $n_0$ and $b$ at level $\lambda$.

We begin with the easiest case, where we act to ensure that every $T^n$ is Souslin. Suppose that $D_\lambda$ codes the triple $\<0,n_0,A>$, where
$n_0$ is a natural number and $A$ is a maximal antichain in $T^{n_0}|\lambda$. In this case, extend the controller trees by specifying any covering
set of branches $B^{n,m}$, and let $G^{n_0}_\lambda$ be the corresponding group of automorphisms of $T^{n_0}|\lambda$ obtained by the compositions
$\pi_{\vec s}$ leading from $n_0$ to $n_0$. By Lemma \ref{lem:basic-sealing-lemma}, we may find a branch $b\in[T^{n_0}|\lambda]$ so that the
resulting set of branches $G^{n_0}_\lambda[b]$ seals $A$. This procedure will ultimately ensure that the tree $T^{n_0}$ is Souslin, because for any
antichain $A\of T^{n_0}$, there would be a stationary set of $\alpha$ such that $D_\alpha$ codes $\<0,n_0,A\intersect T^{n_0}|\alpha>$, at which
point $A$ would be sealed below level $\alpha$. Therefore, the ultimate tree $T^{n_0}$ could have no uncountable antichain.

Similarly, we act to seal antichains in the controller trees $T^{n,m}$. Specifically, if $D_\lambda$ codes $\<1,n,m,A>$ where $A\of T^{n,m}|\lambda$
is a maximal antichain, then we choose the branches $B^{n,m}$ so as to seal $A$, and then choose the other branches $B^{k,r}$, $n_0$ and
$b\in[T^{n_0}|\lambda]$ arbitrarily. The result is that $A$ is sealed in $T^{n,m}$ and so every controller tree $T^{n,m}$ will be Souslin. In this
case, the choice of $n_0$ is irrelevant.

More generally, we act next to ensure that the finite product forcing $T^{0,\ldots,m}=T^{0,1}\cross\cdots\cross T^{m-1,m}$ will have the countable
chain condition. Specifically, suppose that $D_\lambda$ codes $\<2,m,A>$, where $A$ is a maximal antichain in $T^{0,\ldots,m}|\lambda$. We will build
covering sets of branches $B^{i,i+1}\of[T^{i,i+1}|\lambda]$ for $i< m$ to fulfill the requirement that the corresponding products $\vec
b=\<b_0,\ldots,b_{m-1}>$ with $b_i\in B^{i,i+1}$ all lie above a node in $A$. To do this, consider the partial order
\[ \P=(T^{0,1}|\lambda)^\omega\times\cdots\times(T^{m-1,m})^\omega \]
with finite support in every factor. Let us view elements of $\P$ as vectors $\vq=\<q^0,\ldots,q^{m-1}>$, where each $q^i$ is a function from
$\omega$ to $T^{i,i+1}$. We will construct a pseudo-generic filter $G$ in $\P$ by meeting certain dense sets, and then construct the covers
$B_{i,i+1}$, for $i<m$, by setting:
\[ B_{i,i+1}=\{\bigcup\{q^i(k)\st \<q^0,\ldots,q^{m-1}>\in G\}\st
k<\omega\}. \] We ensure that $B_{l,l+1}$ covers $T^{l,l+1}|\lambda$ by meeting the following dense sets, for every $p\in T^{l,l+1}|\lambda$ and
$l<m$.
\[ D_{p,l}=\{\<q^0,\ldots,q^{m-1}>\in\P\st\exists k<\omega\quad p\le
q^l(k)\},\] We ensure that $A$ is sealed by meeting the following dense sets, for every $\vec l\in\omega^m$.
\[ D_{\vec l}=\{\<q^0,\ldots,q^{m-1}>\in\P\st\exists r\in A\quad
\<q^0,\ldots,q^{m-1}>\ge r\}, \] And we ensure that the new branches are cofinal by meeting the following dense sets, for every $i<m$, $j<\omega$ and
$\alpha<\lambda$.
\[ D_{i,j,\gamma}=\{\<q^0,\ldots,q^{m-1}>\in\P\st|q^i(j)|\ge\alpha\},\]

After this, for the other values of $n$ and $k$, we let $B^{n,k}\of[T^{n,k}|\lambda]$ be any covering set of branches and choose $n_0$ and
$b\in[T^{n_0}|\lambda]$ arbitrarily to complete the definition of the $\lambda^{\rm th}$ level as earlier. This part of the construction will ensure
that the ultimate $T^{0,\ldots,m}$ is c.c.c., because any maximal antichain $A\of T^{0,\ldots,m}$ will be anticipated by some $D_\lambda$, and sealed
at stage $\lambda$ as above. So all such antichains are bounded.

Next, we act to ensure that all the trees $T^n$ are rigid. For convenience, we will ensure that they have the unique branch property. Suppose
$D_\lambda$ codes $\<2,n_0,f>$, where $f$ is a potential additional branch for $T^{n_0}|\lambda$. In this case, we extend the controller trees by
specifying any covering set of branches $B^{n,m}$ and let $G^{n_0}_\lambda$ be the corresponding group of automorphisms $\pi_{\vec s}$ of
$T^{n_0}|\lambda$. By Lemma \ref{lem:basic-sealing-lemma}, we may find a branch $b\in[T^{n_0}|\lambda]$ such that corresponding generated set of
branches $G^{n_0}_\lambda[b]$ seals $f$. It follows that $T^{n_0}$ will have the unique branch property (and hence be rigid), since we will have
killed off any potential additional branch.

Similarly, we can also ensure that each controller tree $T^{n,m}$ has the unique branch property. If $D_\lambda$ codes $\<3,n,m,f>$, where $f$ is a
potential additional branch of $T^{n,m}|\lambda$, then we simply choose the branches $B^{n,m}$ so as to seal $f$, and then choose $n_0$ and
$b\in[T^{n_0}]$ arbitrarily. This kills off any potential additional branch for $T^{n,m}$, and so it will have the unique branch property.

Let us now observe some facts that are already determined about the trees we will ultimately construct. Forcing with the controller tree $T^{n,m}$
produces a generic cofinal branch $s\in[T^{n,m}]$. By the inductive assumption 3, the corresponding automorphism $\pi_s$ of ${}^{<\omega_1}2$ takes
$T^n$ to $T^m$. Thus, $\pi_s$ shows that $T^n\cong T^m$ in the forcing extension $V[s]$. So forcing with the $n$-fold product
$T^{0,1,\ldots,m}=T^{0,1}\cross\cdots\cross T^{m-1,m}$ will ensure $T^0\cong T^1\cong\cdots\cong T^m$ in the extension $V[s_0,\ldots,s_n]$. In order
to prove the theorem, therefore, we must ensure that the trees $T^n$ exhibit no other unwanted isomorphism relations. In the next case of our
recursive construction, therefore, we will anticipate and seal any such unwanted potential isomorphisms.

Specifically, suppose that $D_\lambda$ codes $\<5,m,i,j,f>$, where $m\leq i<j$ and $f$ is a $(T^{0,1,\ldots,m}|\lambda)$-potential isomorphism of
$T^{i}|\lambda$ to $T^{j}|\lambda$. That is, $f$ is a function with domain $T^{0,1,\ldots,m}|\lambda$, such that $f(\vec q)$ is a partial isomorphism
of $T^{i}|\lambda$ to $T^{j}|\lambda$, and for any condition $\vec q\in T^{0,1,\ldots,m}|\lambda$, there is a dense set of stronger conditions $\vec
r$, whose $f(\vec r)$ extends $f(\vec q)$ so as to insert any given node of $T^{i}|\lambda$ into the domain and any given node of $T^{j}|\lambda$
into the range. Such functions $f$ exactly arise from $T^{0,1,\ldots,m}$-names for isomorphisms of $T^{i}$ to $T^{j}$. We will now seal this
potential isomorphism. First, by a simple diagonalization meeting countably many dense sets, we choose branches $\vec b=\<b_0,\ldots,b_{m-1}>$ with
$b_k\in[T^{k,k+1}|\lambda]$ in such a way that $f[\vec b]=\bigcup\{f(\vec p)\mid \vec p\leq\vec b\}$ is an isomorphism of $T^i|\lambda$ to
$T^j|\lambda$. According to our construction pattern, we have to specify $n_0$ and a branch $b$ through $T^{n_0}|\lambda$. Set $n_0=i$ and fix an
arbitrary cofinal branch $b\in[T^i|\lambda]$. Next, we will choose further covering sets $B^{n,k}$ of branches for each tree $T^{n,k}|\lambda$ in
such a way so as to guarantee the following crucial property: for any finite list $n_0,\ldots,n_k$ from $n_0=i$ to $n_k=j$ and any $\vec
s=\<s_0,\ldots,s_k>$ with $s_r\in B^{n_r,n_{r+1}}$ we have $\pi_{\vec s}(b)\neq f[\vec b][b]$. If we can do this, then the corresponding definition
as above, of $T^{n,k}(\lambda)=B^{n,k}$ and $T^n(\lambda)$ as the set of all the relevant $\pi_{\vec t}(b)$, will kill off $f[\vec b]$ as an
isomorphism of $T^i$ to $T^j$, since the branch $b$ was extended in $T^i$ but its image $f[\vec b][b]$ was not extended in $T^j$. The key observation
is that we can find the branches $B^{n,k}$ with the crucial property that we mentioned by meeting countably many dense sets. Specifically, we will
associate to each node $p\in T^{n,k}|\lambda$ a branch $b_p\in [T^{n,k}|\lambda]$, determined by our specifying longer and longer initial segments of
it. For each of the countably-many possible patterns $n_0,\ldots,n_k$ leading from $i$ to $j$ and the possible $s_0,\ldots,s_k$ chosen from the $b_p$
we are specifying, we look at the partial information about $\pi_{\vec s}(b)$ that has been determined by the initial segments of $\vec s$. Since the
initial branches $b_0,\ldots, b_m$ that decide $f$ do not involve a controller of $T^j$, it must be at any stage of the construction that at least
one of the branches in such an $\vec s$ is only partially specified. Extending this branch in incompatible ways leads to incompatible isomorphisms
$\pi_{\vec s}$, and so in particular, by making a greater commitment to the partially specified branch $b_p$ in $\vec s$, we can ensure that the
corresponding $\pi_{\vec s}$ acts on $b$ in a way that is different from $f[\vec b][b]$. By enumerating these requirements in an $\omega$-sequence,
we can systematically meet them all, and thereby build the sets of covering branches as desired. It now follows that forcing with $T^{0,\ldots,m}$
will not create unwanted isomorphisms from $T^i$ to $T^j$ when $m\leq i<j$.

Finally, a similar method can be used to ensure that forcing with $T^{0,\ldots,m}$ does not add automorphisms of any $T^i$. Specifically, if there
were a $T^{0,\ldots,m}$-name for a new automorphism of $T^i$, then there would be a $T^{0,\ldots,m}$-potential automorphism of $T^i$, a function $f$
mapping $T^{0,\ldots,m}$ to partial automorphisms of $T^i$, in such a way that one can add nodes to the domain or range of $f(\vec q)$ by
strengthening $\vec q$ to $\vec r$ in $T^{0,\ldots,m}$ and considering $f(\vec r)$. If $D_\lambda$ codes $\<6,m,i,f>$ where $f$ is such a
$(T^{0,\ldots,m}|\lambda)$-potential automorphism of $T^i|\lambda$, then we can seal $f$ as follows: In our construction pattern, again let $n_0=i$.
Fix a branch $b\in[T^i|\lambda]$ and branches $\vb=\<b_0,\ldots,b_{m-1}>$, where $b_i\in[T^{i,i+1}|\lambda]$ in such a way that $f[\vb]$ is a
nontrivial automorphism of $T^i|\lambda$. Now find covering branches $B^{n,k}$ for $T^{n,k}|\lambda$ and a branch as above, by specifying longer and
longer initial segments, so that $f[\vec b][b]\neq\pi_{\vec s}(b)$ for any $\vec s=\<s_0,\ldots,s_m>$ leading from $i$ to $i$. The key point again is
that any list $n_0,\ldots,n_k$ leading from $i$ to $i$ with nontrivial $\vec s$ and $s_r\in B^{n_r,n_{r+1}}$ will necessarily involve at least one
branch $s_r$ that is only partially specified, since one cannot use the initial branches $b_0,\ldots,b_m$ alone to form such a closed loop from $i$
to $i$ without repeating the branches (which would then cancel as the $\pi_s$ commute and are self-inverse).

This completes the recursive construction of the trees $T^n$ and the controller trees $T^{n,m}$. By design, any branch $s$ through $T^{n,m}$ will add
an isomorphism $\pi_s$ from $T^n$ to $T^m$, and so forcing with $T^{0,\ldots,m}$ ensures $T^0\cong\cdots\cong T^m$. We also arranged that this
forcing is essentially a Souslin tree, and hence c.c.c.~and countably distributive, and adds no new automorphisms of any $T^i$ and no unwanted
isomorphisms from $T^i$ to $T^j$ when $m\leq i<j$. So the theorem is proved. \qed

Parts of this construction are redundant, and we don't actually need all the cases that we mentioned. For example, we needn't act explicitly to
ensure that the $T^n$ are rigid, since this will follow by our acting to ensure that they are $T^{0,\ldots,m}$-absolutely rigid. Also, the controller
trees will automatically have the unique branch property, since otherwise they would create unwanted isomorphisms or automorphisms. We included these
redundant simpler cases in the construction simply because they may help to explain the later more complicated parts of the construction, and they
certainly do no harm.

The construction suggests a degree of flexibility that might allow one to improve Theorem \ref{Theorem.DiamondImpliesTrees} by increasing $\omega$ to
any countable ordinal or even $\omega_1$. For example, already the controller trees can be used to make any finitely many trees isomorphic, and one
could arrange that this would not add unwanted isomorphisms or automorphisms. What is needed is a more complicated construction that would anticipate
names for additional automorphisms after forcing with {\it infinite} products of controller trees. This would produce groups whose automorphism tower
can be forced to have any countable ordinal as its height, thereby improving Theorem \ref{Theorem.DiamondImpliesGroups} from $\omega$ to $\omega_1$
as well.

It appears that our methods may generalize to the case of higher cardinals, producing suitably rigid Souslin $\kappa^+$-trees from a suitable
$\diamondsuit_{\kappa^+}$ hypothesis, ultimately constructing groups, whose automorphism towers are highly malleable by forcing. We leave this idea
for a subsequent project.

\bibliographystyle{alpha}

\bibliography{literatur,MathBiblio,HamkinsBiblio}
\end{document}